\documentclass[]{article}

\usepackage{amsfonts}
\usepackage{amscd}
\usepackage{epsfig}
\usepackage{graphicx}
\usepackage[all]{xy}
\usepackage{xypic}
\usepackage{psfrag}
\usepackage{amsmath}
\usepackage{amsthm}
\usepackage{setspace}
\usepackage{url}
\usepackage{here}

\title{Covers of Elliptic Curves and the Lower Bound for Slopes of Effective Divisors on $\overline{\mathcal M}_{g}$}
\author{Dawei Chen}
\date{}

\newtheorem{theorem}{Theorem}[section]
\newtheorem{lemma}[theorem]{Lemma}

\newtheorem{remark}[theorem]{Remark}
\newtheorem{claim}[theorem]{Claim}
\newtheorem{conjecture}[theorem]{Conjecture}
\newtheorem{question}[theorem]{Question}

\begin{document}
\bibliographystyle{plain}
\maketitle

\begin{abstract}
Consider genus $g$ curves that admit degree $d$ covers to elliptic
curves only branched at one point with a fixed ramification type.
The locus of such covers forms a one parameter family $Y$ that
naturally maps into the moduli space of stable genus $g$ curves
$\overline{\mathcal M}_{g}$. We study the geometry of $Y$, and
produce a combinatorial method by which to investigate its slope,
irreducible components, genus and orbifold points. As a by-product
of our approach, we find some equalities from classical number
theory. Moreover, a correspondence between our method and the
viewpoint of square-tiled surfaces is established. We also use our
results to study the lower bound for slopes of effective divisors
on $\overline{\mathcal M}_{g}$.
\end{abstract}

\tableofcontents

\section{Introduction}
Let $\overline{\mathcal M}_{g}$ denote the Deligne-Mumford moduli
space of genus $g$ stable curves. Our aim is to study the
effective cone of $Pic(\overline{\mathcal M}_{g})\otimes \mathbb
Q$. Recall that $Pic(\overline{\mathcal M}_{g})\otimes \mathbb Q$
is generated by the Hodge class $\lambda$ and boundary classes
$\delta_{i}, i = 0,1,\ldots, [{\displaystyle\frac{g}{2}}].$ The
study of effective divisors on $\overline{\mathcal M}_{g}$ has a
long history dating back to \cite{HMu}, \cite{H} and \cite{EH}.
The most important result in this circle of ideas is the
following.

\begin{theorem}
$\overline{\mathcal M}_{g}$ is of general type when $g\geq 24$.
\end{theorem}

The idea of the proof is to write the canonical divisor in the
form $K_{\overline{\mathcal M}_{g}} = A + E$, where $A$ and $E$
are ample and effective $\mathbb Q$-divisors respectively. Then
$K$ is big and $\overline{\mathcal M}_{g}$ is of general type. In
practice, we know that $K_{\overline{\mathcal M}_{g}} = 13\lambda
- 2\delta - \delta_{1}$ and $\lambda - \alpha\delta$ is ample for
$0 <\alpha \ll 1$, where $\delta$ is the total boundary. So as
long as we obtain an effective divisor $D = a\lambda -
{\displaystyle\sum_{i=0}^{[\frac{g}{2}]}b_{i}\delta_{i}}, a, b_{i}
> 0$ with ${\displaystyle\frac{a}{b_{i}}} <
{\displaystyle\frac{13}{2}}$ for all $i$ and
${\displaystyle\frac{a}{b_{1}}}<{\displaystyle\frac{13}{3}}$, then
we can write $K_{\overline{\mathcal M}_{g}} $ in the desired form
$A+E$. Define the number
 ${\displaystyle\frac{a}{min \ \{b_{i}\}}}$ to be the slope $s(D)$ of the divisor $D$. Then the remaining task is to
 find effective divisors with small slope. \\

 In \cite{EH}, the famous Brill-Noether divisor was shown to do the trick. Recall that when the Brill-Noether number $\rho(g,r,d)=g-(r+1)(g-d+r)$ is equal to $-1$, the Brill-Noether divisor $BN^{r}_{d}$ consists of genus $g$ curves which possess a $g^{r}_{d}$. Amazingly, $s(BN^{r}_{d}) = 6+{\displaystyle\frac{12}{g+1}}$, the so-called Brill-Noether bound, is independent of $d, r$, and is less than ${\displaystyle\frac{13}{2}}$ if $g\geq 24$. For almost
 20 years after the publication of \cite{EH}, no one succeeded in finding any effective divisor of slope lower than the Brill-Noether bound.
This led to the formulation, in \cite{HM}, of the following slope
conjecture.

 \begin{conjecture}
The  Brill-Noether bound $6+{\displaystyle\frac{12}{g+1}}$ provides a lower bound for slopes of effective divisors on $\overline{\mathcal M}_{g}$. The bound is sharp if and only if $g+1$ is composite.
 \end{conjecture}

Unfortunately, the slope conjecture is false. In \cite{FP}, Farkas
and Popa study genus 10 curves contained in $K3$ surfaces, or
equivalently curves embedded in $\mathbb P^{4}$ by $g^{4}_{12}$'s
whose images are contained in quadric threefolds. These curves
sweep out a divisor $\mathcal K$ in $\overline{\mathcal M}_{10}$.
Surprisingly, $s(\mathcal K) = 7$, which is less than the
Brill-Noether bound ${\displaystyle 7\frac{1}{11}}$. So $\mathcal
K$ serves as a first counterexample to the slope conjecture.
Recently, in \cite{Kh}, \cite{F1},
\cite{F2} and \cite{F3}, more counterexamples are constructed, and in these papers Farkas also announces that $ \overline{\mathcal M}_{22}$ and $\overline{\mathcal M}_{23}$ are of general type. \\

The above results only tell us one side of the story, where slopes
of effective divisors are concerned. On the other side, to the
best of the author's knowledge, all effective divisors that have
been studied thus far have slope greater than 6. So in some sense,
the slope conjecture is not far from the truth. A natural question
is to ask:

\begin{question}
Is there a lower bound for slopes of effective divisors on
$\overline{\mathcal M}_{g}$?
\end{question}

A uniform lower bound $\varepsilon$ for all $g$ would enable us to solve the Schottky problem, which roughly says that the modular forms of slope $\geq \varepsilon$ cut out $\mathcal M_{g}$ in $\mathcal A_{g}$, the moduli space for abelian varieties of dimension $g$, see, e.g., \cite{F4}. \\

The idea to obtain a lower bound is to construct moving curves in
$\overline{\mathcal M}_{g}$. An irreducible curve $C\in
\overline{\mathcal M}_{g}$ is called a moving curve if the family
of its deformations fills in a Zariski open set of
$\overline{\mathcal M}_{g}$. We can define the slope of $C$ by
$s(C) = {\displaystyle\frac{C\ldotp\delta}{C\ldotp\lambda}}$. If
an effective divisor $D$ does not contain $C$, then $D\ldotp C\geq
0$ and $s(D)\geq s(C).$ Hence, the smallest slope of curves in
this family bounds the slope of effective divisors, since a
divisor cannot contain all the deformations of a moving curve!
Along these lines, in \cite{HM}, Harris and Morrison study genus
$g$ curves as degree $d$ covers of $\mathbb P^{1}$ for large $d$.
If we vary the branch points in $\mathbb P^{1}$, we just obtain
moving curves in $\overline{\mathcal M}_{g}$. The slopes of such
moving curves do provide good lower bounds when $g$ is small.
However, as $g\rightarrow \infty$, these bounds tend to 0.
Therefore, the existence of a uniform positive lower bound for the slope for all $g$ remains unknown. \\

The same circle of ideas is mentioned in \cite{CHS}. Instead of
using covers of $\mathbb P^{1}$, Coskun, Harris and Starr choose
canonical curves in $\mathbb P^{g-1}$ to study. In order to get
moving curves in $\overline{\mathcal M}_{g}$, they impose on
canonical curves some conditions like contacting fixed linear
subspaces with suitable orders. This method works pretty well for
genus up to 6. But the general case is hard to analyze since we do
not have good knowledge about the geometry of canonical curves and
the Gromov-Witten theory for high genus curves is far from complete. \\

Here we push the idea in \cite{HM} further to consider covers of
elliptic curves. Take a cover $\pi: C\rightarrow E$ that only has
one branch point at the marked point $p$ of the elliptic curve
$E$. If we vary the complex structure of $E$, $C$ also varies and
its locus in $\overline{\mathcal M}_{g}$ forms a 1-dimensional
subscheme $Y$. The geometry of $Y$ is our main interest.
We need to do three things as follows. \\

Firstly, we want to calculate the slope of $Y$. Here we use a very
classical approach originally due to Hurwitz. He used 
symmetric groups to enumerate covers of $\mathbb P^{1}.$
Similarly, the information of a degree $d$ branched cover of an
elliptic curve can be encoded in a certain equivalence class of
solution pairs in $S_{d}$. We fix the ramification type of the
cover by a conjugacy class $\sigma = (l_{1})\cdots (l_{m})$ in
$S_{d}$, i.e.,
$\pi^{-1}(p)={\displaystyle\sum_{i=1}^{m}l_{i}q_{i}}$. Let
$(1^{a_{1}}2^{a_{2}}\dots d^{a_{d}})$ also denote a conjugacy
class of $S_{d}$ whose number of length $i$ cycles equals $a_{i}$,
$\sum\limits^{d}_{i=1} ia_{i} = d$. For two pairs $(\alpha,
\beta), (\alpha', \beta') \in S_{d}\times S_{d}$, we define an
equivalence relation by $(\alpha, \beta) \sim (\alpha', \beta')$
if there exists an element $ \tau \in S_{d}$ such that $\tau
\alpha \tau^{-1} = \alpha'$ and $\tau \beta \tau^{-1} = \beta'$.
Moreover, let $<\alpha, \beta>$ denote the subgroup of $S_{d}$
generated by $\alpha$ and $\beta$. Define a set
$$Cov^{g,d,\sigma} = \{(\alpha, \beta)\in S_{d} \times S_{d}\ |\ \alpha\beta\alpha^{-1}\beta^{-1}\in \sigma, $$  $$ <\alpha, \beta>
\mbox{is a transitive subgroup of} \ S_{d}\}, $$ and its subset
$$Cov_{1^{a_{1}} 2^{a_{2}} \dots d^{a_{d}}}^{g,d,\sigma} = \{(\alpha, \beta)\in Cov^{g,d,\sigma}\ |\ 
\beta \in (1^{a_{1}} 2^{a_{2}} \dots d^{a_{d}})\}$$ based on the
conjugacy type of $\beta$. Then the cover $\pi$ corresponds to one
equivalence class of solution pairs in $Cov^{g,d,\sigma} $. Readers should bear in mind that for this cover, 
$\alpha, \beta$ correspond to the monodromy images of a standard symplectic basis of $H_{1}(E)$. Pick a loop $\gamma
=\alpha\beta\alpha^{-1}\beta^{-1}$ around the branched point $q$. The monodromy image of $\gamma$ has to be an element 
in the ramification class $\sigma$. The transitivity condition imposed on $<\alpha, \beta>$ is because we want the cover to be connected. See the following picture. \\

\begin{figure}[H]
    \centering
    \psfrag{a}{$\alpha$}
    \psfrag{b}{$\beta$}
    \psfrag{E}{$E$}
    \psfrag{q}{$q$}
    \psfrag{r}{$\gamma$}
    \includegraphics[scale=0.4]{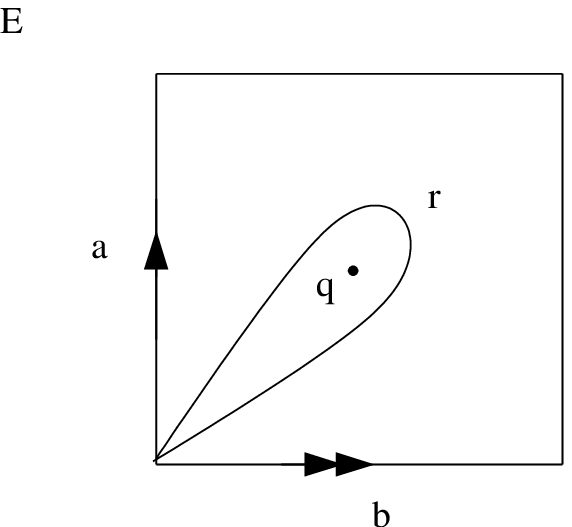}
\end{figure}

Let $N_{1^{a_{1}} 2^{a_{2}}\dots d^{a_{d}}}$ denote the order of the
set $Cov^{g,d,\sigma}_{1^{a_{1}} 2^{a_{2}} \dots d^{a_{d}}}$
modulo the equivalence relation, i.e.,
$$N_{1^{a_{1}} 2^{a_{2}}\dots d^{a_{d}}} = |Cov_{1^{a_{1}} 2^{a_{2}} \dots d^{a_{d}}}^{g,d,\sigma}/\sim|.$$
When there is no ambiguity, we simply use $Cov$ to denote $Cov^{g,d,\sigma}$.
Let us further define a total sum $$N = \sum_{1^{a_{1}}
2^{a_{2}}\dots d^{a_{d}}} N_{1^{a_{1}} 2^{a_{2}}\dots
d^{a_{d}}}$$ and a weighted sum $$M =
\sum_{1^{a_{1}} 2^{a_{2}}\dots d^{a_{d}}}
(\frac{a_{1}}{1}+\frac{a_{2}}{2}+\dots +
\frac{a_{d}}{d})N_{1^{a_{1}} 2^{a_{2}}\dots d^{a_{d}}}.$$
It is clear that $N$ and $M$ only depend on $g, d$ and $\sigma$.
Now we can state the main result about the slope of $Y$.
\begin{theorem}
\label{slope}
$$s(Y)= \frac{12M}{M+\left(d-{\displaystyle\sum_{i=1}^{m}\frac{1}{l_{i}}}\right){\displaystyle\frac{N}{12}}}$$
\end{theorem}

\begin{remark}
Note that $s(Y)$ only depends on the ratio
${\displaystyle\frac{N}{M}}$.
\end{remark}

Secondly, we have to study the irreducible components of $Y$. By
using the symmetric group $S_{d}$ again, we can see that two
covers lie in the same component of $Y$ if and only if the
solution pair corresponding to one can be sent to the solution
pair corresponding to the other via a monodromy action. More
precisely, define two actions $a:(\alpha,\beta)\rightarrow
(\alpha, \alpha\beta)$ and $b:(\alpha, \beta)\rightarrow
(\alpha\beta, \beta)$ acting on the set $Cov$. It is not hard to
check that both actions are well defined for $Cov/\sim$. Let $G$
denote the group of actions generated by $a, b$ and call $G$ the
monodromy group of $Y$. We give the following criterion.
\begin{theorem}
\label{monodromy} Take two covers represented by two solution
pairs $(\alpha,\beta)$ and $(\alpha',\beta')$ in $Cov/\sim$. Then
the points in $Y$ corresponding to these two covers are contained in
the same component of $Y$ if and only if there exists an element
in the monodromy group $G$ that sends $(\alpha,\beta)$ to
$(\alpha',\beta')$.
\end{theorem}

Finally, we have to verify whether or not $Y$ is a moving
curve in $\overline{\mathcal M}_{g}$. Readers may already be aware
of the fact that $Y$ is more likely rigid, since a cover
corresponding to a point in $Y$ only has one branched point, whose
variation along the target elliptic curve does not change the
moduli of this elliptic curve. Nevertheless, technically we do not
need moving curves. Instead, when the degree $d$ of the
covers varies, as long as the union of all $Y_{g,d,\sigma}$ is
Zariski dense in $\overline{\mathcal M}_{g}$, an effective divisor
cannot contain all of them. So the smallest slope
$s(Y)$ in this union still provides a lower bound for
slopes of effective divisors. Luckily, for ramification class
$\sigma$ of suitable types, we have an affirmative answer to this
density problem.
\begin{theorem}
\label{density} When $g$ and $\sigma = (l_{1})\cdots (l_{m})$ are
fixed but $d$ varies, the image of the union of $Y_{g,d,\sigma}$ is dense in
$\overline{\mathcal M}_{g}$ if and only if the number of
ramification points is less than $g$, i.e., $|\{l_{i}: l_{i}\geq
2\}| < g$.
\end{theorem}

We apply the above results to the case $g=2$. About the slope, we
have the following conclusion.
\begin{theorem}
\label{g2s} When $g = 2$, $\sigma$ can be either of type
$(3^{1}1^{d-3})$ or of type $(2^{2}1^{d-4})$. In both cases $Y$
has constant slope $s(Y)= 10$ independent of $d$.
\end{theorem}

Note that 10 is the sharp lower bound for slopes of effective
divisors on $\overline{\mathcal M}_{2}$, cf. Remark \ref{Mc}. We
will give a direct proof of Theorem \ref{g2s} in section 4 but
only for prime number $d$. The general case follows from an
indirect argument in Remark \ref{Mc}.\\

We can further calculate the genus of $Y$. The following results
are only for prime $d$.
\begin{theorem}
\label{g21m} When $g = 2$, $d$ is prime and $\sigma$ is of type
$(3^{1}1^{d-3})$,
$$g(Y) =  1+\frac{1}{8}(d-1)(d-2)(15d+23) - 6\Big(\sum_{a_{1}l_{1}+a_{2}l_{2}=d\atop l_{1}>l_{2}}(l_{1}, a_{1})
(l_{2}, a_{2})\Big), $$ where $(m,n)$ denotes the greatest common
divisor of $m$ and $n$. In particular, asymptotically $g(Y)\sim
{\displaystyle\frac{15}{8}d^{3}}$.
\end{theorem}

\begin{theorem}
\label{g22m} When $g = 2$, $d$ is prime and $\sigma$ is of type
$(2^{2}1^{d-4})$,
$$g(Y) =  1+\frac{1}{12}(d-1)(d-3)(10d^{2}-13d-14)-6\Big(\sum_{a_{1}l_{1}+a_{2}l_{2}+a_{3}l_{3}=d\atop l_{1}=l_{2}+l_{3}>l_{2}>l_{3}}\prod_{i=1}^{3}(l_{i}, a_{i}) $$
$$-\sum_{a_{1}l_{1}+a_{2}l_{2}=d\atop l_{1}>l_{2}}(l_{1}-2)(l_{1}, a_{1})(l_{2}, a_{2}) - \sum_{2a_{2}+a_{1}=d}\frac{a_{1}-1}{(a_{2},2)}\Big).  $$
In particular, asymptotically $g(Y)\sim {\displaystyle\frac{5}{6}d^{4}}$.
\end{theorem}

Theorems \ref{g21m} and \ref{g22m} will be proved in section 5 by considering local monodromy actions. \\

When $g$ is 3, the combinatorics becomes harder. However, we check
one case by computers: when $\sigma$ is of type $(5^{1}1^{d-5})$
and $d$ is prime, it looks like that $s(Y)$ tends to 9, cf.
Conjecture \ref{g3s}. Note that 9 is the sharp lower bound for
slopes of effective divisors on $\overline{\mathcal M}_{3}$, cf. \cite{HM}. \\

Unfortunately, although in principle by Theorems \ref{slope} and \ref{monodromy} we know almost everything about $Y$,
in practice it seems that the combinatorial problems we encounter are quite complicated, especially when $g$ is large. It would be very interesting to obtain
more numerical evidences to see if this method does provide a lower bound for slopes of effective divisors.   \\

This paper is organized in the following way. In section 2 we
study the geometry of $Y$, namely we consider slope, monodromy and density.
Then we apply the results to the genus 2 case and several examples of the genus 3 case in section 3. In section 4,
a general method is introduced to express $N$ and $M$ in terms of some generating functions.
Then we study the genus and orbifold points of $Y$ in section 5.
Afterwards, a correspondence between our viewpoint and the viewpoint of square-tiled surfaces is established in section 6. Finally in the appendix,
we prove some equalities of certain series appearing in section 5 by a trick of classical number theory.   \\

{\bf Acknowledgements.} I am grateful to my advisor Joe Harris who
told me about this problem and gave me many valuable suggestions.
I would also like to thank Sabin Cautis, Izzet Coskun, Noam
Elkies, Thomas Lam and Curtis T. McMullen for useful
conversations.

\section{The Global Geometry of $Y$}
In this section, we will focus on the slope, monodromy and density
of $Y$ and prove Theorems \ref{slope}, \ref{monodromy} and
\ref{density} one by one. First, let us give a rigorous construction of $Y$. \\

Let $X \cong \mathbb P^{1}$ parameterize a general pencil of plane
cubics. Blow up the 9 base points to obtain a smooth surface $S$
that is an elliptic fibration over $X$. There are 12 rational
nodal curves as special fibers over $b_{1},\ldots,b_{12}\in X$.
Fix a section $\Gamma$ corresponding to the blow up of one of the
9 base points. With this section, $X$ can be considered as a
12-sheeted cover of the moduli space of elliptic curves
$\overline{\mathcal M}_{1,1}$. Take a general fiber $(E,p)$ over
$b\in X$. Consider all the possible degree $d$ covers $\pi:
C\rightarrow E$ from a genus $g$ connected curves $C$ to $E$
branched only at $p$. Let $\pi^{-1}(p)= l_{1}q_{1}+ \ldots
+l_{m}q_{m}$, where $(l_{1})\cdots (l_{m})$ is a fixed conjugacy
class $\sigma$ of the symmetric group $S_{d}$. When $b$ varies in
$X_{0}= X\backslash \{b_{1},\ldots,b_{12}\}$, the locus of such
covers also varies in a 1-dimensional Hurwitz scheme
$Y^{0}_{g,d,\sigma}$. Now, if $b$ approaches some point $b_{i}$,
the cover degenerates to a cover of a rational nodal curve, in the
sense of admissible covers. Hence, we can compactify $Y^{0}_{g,d,\sigma}$ by the space of
admissible covers ${Y}_{g,d,\sigma}$. When
there is no ambiguity, we will simply use $Y^{0}$ and $Y$ instead. \\

The curve $Y$ is an $N$-sheeted cover of $X$, possibly branched at
$b_{1},\ldots, b_{12}$. $N$ is equal to the number of distinct
genus $g$ degree $d$ covers of a general plane cubic only branched
at one point of ramification type $\sigma$. Note that this number
$N$ equals the total sum $$N = \sum_{1^{a_{1}} 2^{a_{2}}\dots
d^{a_{d}}} N_{1^{a_{1}} 2^{a_{2}}\dots d^{a_{d}}}$$ defined in the
introduction section.

\subsection{Slope}
There is a natural map from the space of admissible covers $Y$ to
$\overline{\mathcal M}_{g}$. So it makes sense to talk about the
intersections $Y\ldotp \delta$ and $Y\ldotp \lambda$ by pulling
back the corresponding classes to $Y$. Now we can prove
Theorem~\ref{slope}.
\begin{proof} We have to figure out the two numbers $Y\ldotp
\delta$ and $Y\ldotp \lambda$. The former can be worked out by the
following argument. We want to establish a diagram of maps:
$$\xymatrix{
T \ar[r]^\varphi \ar[d]  & S \ar[d] \\
Y \ar[r]^\psi     & X }$$ such that over a point $y \in Y$ mapping
to $x \in X$, the fiber over $y$ is the corresponding cover of the
fiber over $x$, i.e., T is the universal covering curve over $Y$. See the picture below. 

\begin{figure}[H]
    \centering
    \psfrag{x}{$x$}
    \psfrag{y}{$y$}
    \psfrag{X}{$X$}
    \psfrag{Y}{$Y$}
    \psfrag{S}{$S$}
    \psfrag{T}{$T$}
    \psfrag{bi}{$b_{i}$}
    \psfrag{N}{$N$}
    \includegraphics[scale=0.5]{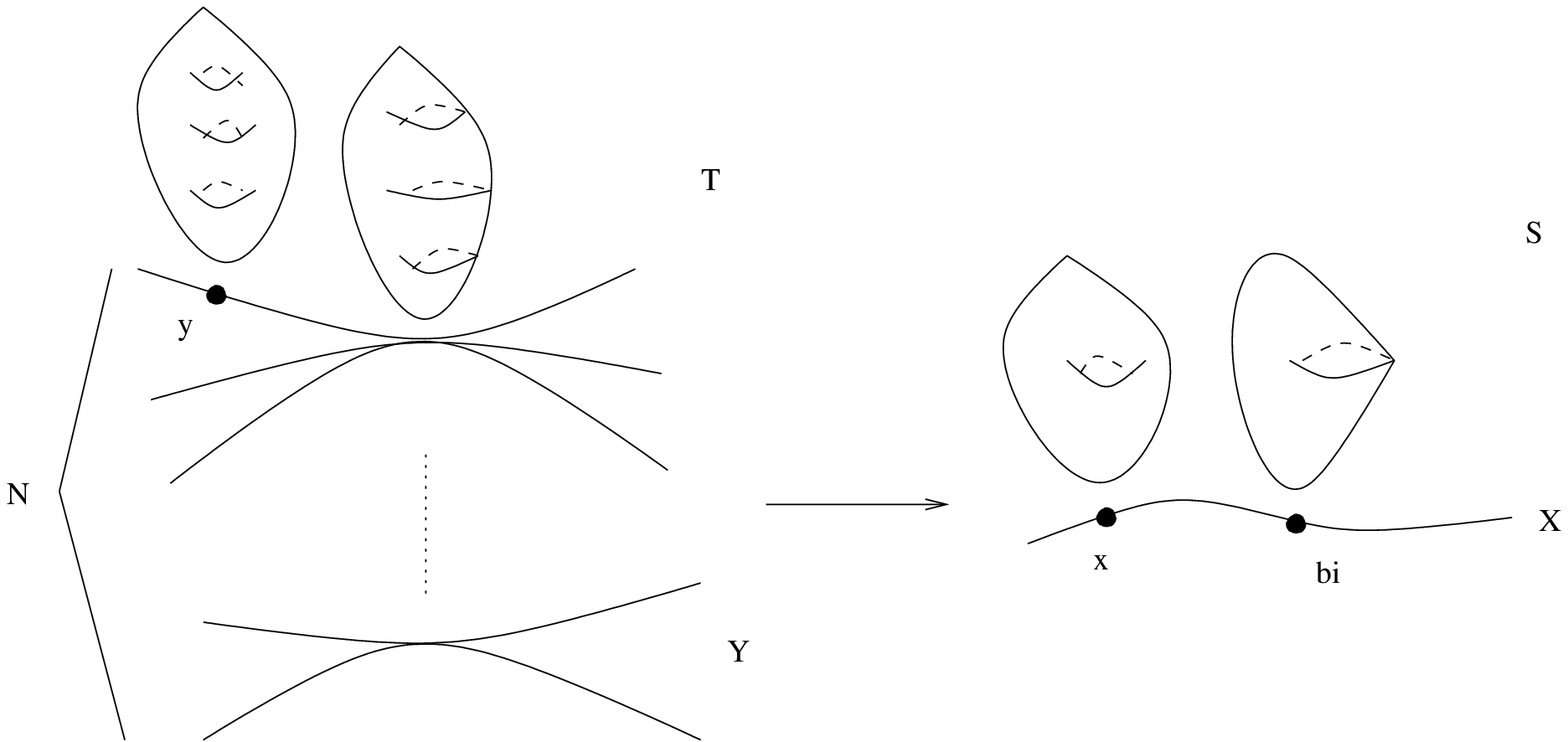}
\end{figure}

When $x$
is in $X^{0}$, locally we can always construct $T$. Now in a small
neighborhood of $b_{i}$, we can identify all the smooth fibers to
a fixed elliptic curve $(E,p)$. Let $\alpha, \beta$ be a standard symplectic 
basis of $\pi_{1}(E,p)$ such that $\beta$ is the vanishing cycle
when the smooth elliptic curve degenerates to the rational nodal
curve over $b_{i}$. We abuse notation and denote their monodromy
images in $S_{d}$ also by $\alpha, \beta$, and assume $\beta \in
(1^{a_{1}}2^{a_{2}}\ldots d^{a_{d}})$. Now the key observation is
that when the elliptic curve degenerates, the combinatorial type of the solution pair corresponding
to the smooth cover determines the type of the degenerate cover.
More precisely, the degenerate covering curve will have $a_{i}$
nodes at which the map $\varphi$ is given locally by $u\rightarrow
x=u^{i}, v\rightarrow y=v^{i}$, where $(u,v)$ and $(x,y)$ are the
local charts of $T$ and $S$ respectively. Each such node
contributes ${\displaystyle\frac{1}{i}}$ to the intersection
$Y\ldotp \delta$. Therefore, we have
$$ Y\ldotp \delta = 12\sum_{1^{a_{1}}
2^{a_{2}}\dots d^{a_{d}}}(\frac{a_{1}}{1}+\frac{a_{2}}{2}+\dots +
\frac{a_{d}}{d}) N_{1^{a_{1}} 2^{a_{2}}\dots d^{a_{d}}}=
12M.$$
The constant 12 comes from the number of singular fibers of $S$. \\

To calculate $Y\ldotp \lambda$, by Mumford's formula $\lambda =
{\displaystyle \frac{\delta + \kappa}{12}}$, it suffices to work
out the intersection $Y\ldotp \kappa = \omega^{2}_{T/Y}$. The
ramification class $\sigma$ is of type $(l_{1})\ldots (l_{m})$ so
we have $m$ sections $\Gamma_{1},\ldots,\Gamma_{m}$ of $T$ such
that $\varphi^{*}\Gamma =
{\displaystyle\sum_{i=1}^{m}l_{i}\Gamma_{i}}$. By Riemann-Hurwitz,
the relative dualizing sheaves $\omega_{T/Y}$ and $\omega_{S/X}$
satisfy the relation
$$\omega_{T/Y}=\varphi^{*}\omega_{S/X}+\sum_{i=1}^{m}(l_{i}-1)\Gamma_{i}.$$
Moreover, $$\varphi_{*}\Gamma_{i}=N\cdotp \Gamma,\
(\varphi^{*}\omega_{S/X})^{2}=dN\cdotp \omega^{2}_{S/X}=0,$$ and
$$\Gamma_{i}\ldotp\Gamma_{j}=0,\  i\neq j.$$ So we get
$$l_{i}\Gamma^{2}_{i} = \Gamma_{i}\ldotp
(\varphi^{*}\Gamma)=(\varphi_{*}\Gamma_{i})\ldotp\Gamma=N(\Gamma^{2})=-N,$$
and $$\Gamma_{i}\ldotp (\varphi^{*}\omega_{S/X})=
(\varphi_{*}\Gamma_{i})\ldotp \omega_{S/X} =
N(\Gamma\ldotp\omega_{S/X})=N(-\Gamma^{2})=N.$$ Now it is routine
to check that $$\omega^{2}_{T/Y} =
\left(\sum_{i=1}^{m}l_{i}-\sum_{i=1}^{m}\frac{1}{l_{i}}\right)
N = \left(d-\sum_{i=1}^{m}\frac{1}{l_{i}}\right) N.$$ So we
get the desired formula for $s(Y)$.
\end{proof}

\begin{remark}
\label{bc} To complete the commutative diagram in the above proof,
when the map is given locally by $u\rightarrow x=u^{i},
v\rightarrow y=v^{i}$ where $(u,v)$ and $(x,y)$ are the local
charts of $S$ and $T$ respectively, a base change of degree
divisible by $i$ is necessary. So we can make a degree $d!$ base
change once for all to realize such maps globally. After the base
change, $Y\ldotp\delta$ and $Y\ldotp \lambda$ both have to be
multiplied by $d!$, so the quotient $s(Y)$ remains the same. \\

A fancier interpretation is that by using a minimal base change,
we pass from the space of admissible covers $Y$ to the stack of
admissible covers $\widetilde{Y}$ where a universal covering map
lives. This viewpoint is crucial when we study the local geometry
of $Y$ in section 5.
\end{remark}

\begin{remark}
\label{in} Recall that in our original definition of a cover and
its automorphism, we did not allow variation of the target
elliptic curve. So an automorphism of a cover $\pi: C\rightarrow
E$ is given by the following commutative diagram,
$$\xymatrix{
C \ar[rr]^\phi \ar[dr]_\pi & & C \ar[dl]^\pi \\
& E  &  }$$
where $\phi$ is an automorphism of $C$. \\

However, the elliptic curve $E$ always has an automorphism induced
by its involution $\iota$ that sends $(\alpha,\beta)$ to
$(\alpha^{-1},\beta^{-1})$. If we further identify the two
solution pairs $(\alpha,\beta)\sim (\alpha^{-1},\beta^{-1})$, that
is, we allow the following commutative diagram to be considered as
an automorphism of the cover $\pi$,
$$\xymatrix{
C \ar[rd]^{\iota\circ\pi} \ar[d]_\pi  &  \\
E \ar[r]_\iota     & E } $$ we would get a new space $Y'$. $Y$ can
be viewed as a double cover of $Y'$. More precisely, for a
component $Z$ of $Y$, if a general cover in $Z$ corresponding to a
solution pair $(\alpha,\beta)$ does not have the automorphism
induced by $\iota$, i.e., there does not exist an element $\tau\in
S_{d}$ such that
$(\tau\alpha\tau^{-1},\tau\beta\tau^{-1})=(\alpha^{-1},\beta^{-1})$,
then $Z$ is a double cover of the corresponding component $Z'$ of
$Y'$. On the other hand, if there is some $\tau\in S_{d}$ such
that
$(\tau\alpha\tau^{-1},\tau\beta\tau^{-1})=(\alpha^{-1},\beta^{-1})$,
then every cover in $Z$ has the automorphism induced by $\iota$.
Hence, $Z$ is a double curve and its reduced structure is the same
as $Z'$. Note for the purpose of slope calculation that $s(Z)$ and
$s(Z')$ are always the same.
\end{remark}

\subsection{Monodromy}
The curve $Y$ may be reducible. Actually if two pairs
$(\alpha,\beta)$ and $(\alpha',\beta')$ generate non-conjugate
subgroups of $S_{d}$, then the corresponding two covers must
be contained in different components of $Y$. \\

Now let us consider the monodromy of the map $Y\rightarrow X$.
More precisely, we want to study the $\pi_{1}$-monodromy map
$\rho_{\pi}: \pi_{1}(X^{0}, b) \rightarrow
Out^{+}(\pi_{1}(X_{b}))= Out^{+}(\mathbb Z*\mathbb Z)$, where $b$
is a fixed base point of $X^{0}$, $X_{b}$ is the fiber over $b$
with one marked point and $Out^{+}(\mathbb Z*\mathbb
Z)=Aut(\mathbb Z*\mathbb Z)/Inn(\mathbb Z*\mathbb Z)$ is the
orientable outer automorphism group.
\begin{lemma}
The map $\rho_{\pi}$ is surjective and its image $Out^{+}(\mathbb
Z*\mathbb Z) \simeq SL_{2}(\mathbb Z)$. Moreover, the group $G$ of
$\pi_{1}$-monodromy acting on $Cov/\sim$ can be generated by two
operations: $a: (\alpha, \beta)\rightarrow (\alpha, \alpha\beta)$
and $b: (\alpha, \beta)\rightarrow (\alpha\beta, \beta)$.
\end{lemma}

\begin{proof}
We give an indirect proof using the fact that the
$H_{1}$-monodromy map $\rho_{H}$ of a general pencil of plane
cubics is
surjective. \\

For $H_{1}$-monodromy, the marked point does not affect homology.
So we consider the map $\rho_{H}: \pi_{1}(X^{0}, b) \rightarrow
Out^{+}(H_{1}(X_{b}))= Aut^{+}(\mathbb Z\oplus\mathbb Z)\simeq
SL_{2}(\mathbb Z)$. In order to show that $\rho_{H}$ is surjective
for a general pencil $X$, it suffices to exhibit a special pencil
for which the claim holds. Actually some examples are explicitly
worked out in \cite{Sa} and the $H_{1}$-monodromy maps are
surjective as we expect. \\

Now for $\pi_{1}$-monodromy, we can construct the following
commutative diagram:
$$\xymatrix{
                                                   & \Gamma_{1,1}\ar[d]_\simeq \ar@/^5pc/[ddd]_\simeq           \\
\pi_{1}(X^{0},b)\ar[r]^{\rho_{\pi}}\ar[ur]\ar[d]_\simeq & Out^{+}(\mathbb Z*\mathbb Z)\ar[d]_{\phi} \\
\pi_{1}(X^{0},b)\ar[r]^{\rho_{H}}\ar[dr]             & Aut^{+}(\mathbb Z\oplus\mathbb Z)         \\
                                                   & \Gamma_{1}\ar[u]^\simeq
                                                   }$$
Here $\Gamma_{1}$ and $\Gamma_{1,1}$ are the mapping class groups
for an ordinary torus and a torus with one marked point
respectively. It is well known that $\Gamma_{1}$ and
$Out^{+}(\mathbb Z*\mathbb Z)$ are both isomorphic to
$SL_{2}(\mathbb Z)$. The map $\phi$ is induced by quotienting out
commutators, i.e., $\mathbb Z*\mathbb Z / ([\alpha, \beta])\simeq
\mathbb Z\oplus\mathbb Z$, where $\alpha$ and $\beta$ are the
generators of the free group $\mathbb Z*\mathbb Z$. Moreover, we
have the isomorphism
$\Gamma_{1,1}\stackrel{\simeq}{\rightarrow}\Gamma_{1}$, cf.
\cite{Bi}.
So the conclusion that $\rho_{\pi}$ is surjective follows from the fact that $\rho_{H}$ is surjective. \\

Finally, $Out^{+}(\mathbb Z*\mathbb Z)$ is generated by the two
basic actions induced from Dehn twists along the two loops
represented by $\alpha$ and $\beta$. So correspondingly the
monodromy actions send $(\alpha, \beta)$ to $(\alpha,
\alpha\beta)$ and $(\alpha\beta, \beta)$ respectively.
\end{proof}

Now Theorem~\ref{monodromy} follows easily from the above lemma.

\begin{remark}
Since $(\alpha, \alpha\beta)\sim \alpha^{-1}(\alpha, \alpha\beta)\alpha = (\alpha, \beta\alpha)$, the
monodromy actions are well defined up to the equivalence relation.
The two actions $a$ and $b$ correspond to $\left(
\begin{array}{cc}
1 & 1 \\
0 & 1 \end{array} \right)$ and $\left( \begin{array}{cc}
1 & 0 \\
1 & 1 \end{array} \right)$ respectively, which generate
$SL_{2}(\mathbb Z)$.
\end{remark}

\begin{remark}
\label{lm}
Assume that $\beta$ corresponds to the local vanishing cycle around a nodal fiber over 
$b_{i} \in X$. By the Picard-Lefschetz formula or Kodaira's classification of elliptic fibrations, going along a small loop around $b_{i}$ once will send a solution pair $(\alpha, \beta)$ to $(\alpha\beta, \beta)$. See the following picture. 

\begin{figure}[H]
    \centering
    \psfrag{X}{$X$}
    \psfrag{b}{$\beta$}
    \psfrag{bi}{$b_{i}$}
    \includegraphics[scale=0.5]{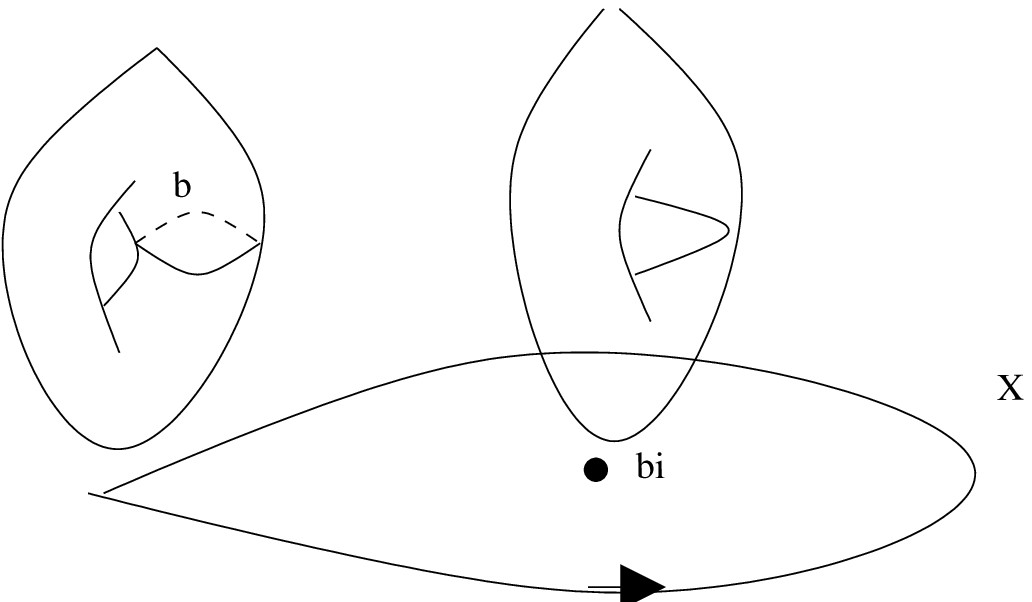}
\end{figure}

We call this action as the local monodromy action, in order to distinguish it from the global monodromy. 
It is useful when we study the genus of $Y$ by applying the Riemann-Hurwitz formula to the map $Y\rightarrow X$, 
cf. section 5. 
\end{remark}

\subsection{Density}
In this section, whenever we mention density, we mean Zariski
density. Notation is the same as before. Let $\sigma$ be the
conjugacy class $(l_{1})\cdots (l_{m})$ of $S_{d}$. We assume that
$l_{1},\ldots,l_{k}$ are greater than 1 and that the rest of the
$l_{i}$'s are all equal to 1. By Riemann-Hurwitz, $2g-2 =
{\displaystyle\sum_{i=1}^{k}(l_{i}-1)= \Big(\sum_{i=1}^{k}l_{i}\Big)-k}$.
Now the pullback of a holomorphic 1-form from the target elliptic
curve becomes a holomorphic 1-form $\omega$ on the covering curve
$C$ and $(\omega)={\displaystyle\sum_{i=1}^{k}\mu_{i}q_{i}}$,
where $\mu_{i}=l_{i}-1$ and
${\displaystyle\sum_{i=1}^{k}\mu_{i}=2g-2}$. So we can fix the
ramification type by fixing $l_{1},\ldots,l_{k}$ but adding more
length-1 cycles to $\sigma$. Hence, the degree $d$
of the map can vary from $2g-2+k$ to infinity. \\

Now we consider the union of all possible covers in the above
sense. An equivalent statement of Theorem~\ref{density} is the
following.
\begin{theorem}
\label{D} 
The image of ${\displaystyle\bigcup_{d=2g-2+k}^{\infty}Y^{0}_{g,d,\sigma}}$ is dense
in $\mathcal M_{g}$ if and only if the inequality 
${\displaystyle\sum_{i=1}^{k}(\mu_{i}-1)\leq g-1}$ holds, i.e., $k\geq
g-1$.
\end{theorem}

\begin{proof}
Let $\mu = (\mu_{1},\ldots,\mu_{k})$ be a partition of $2g-2$. We
consider the moduli space $\mathcal H(\mu)$ parameterizing
$(C,\omega, q_{1},\ldots,q_{k})$ where $[C]\in \mathcal M_{g}$,
$\omega\in H^{0}(K_{C})$, and
${\displaystyle(\omega)=\sum_{i=1}^{k}\mu_{i}q_{i}}$. The
dimension of $H(\mu)$ is
$$3g-3+g-\left(\sum_{i=1}^{k}(\mu_{i}-1)\right)=2g-1+k = n.$$ The space
$H(\mu)$ can locally be equipped with a coordinate system. Pick a
basis $\gamma_{1},\ldots,\gamma_{n}\in
H_{1}(C,q_{1},\ldots,q_{k};\mathbb Z)$ the relative homology of
$C$ with $k$ marked points, such that
$\gamma_{1},\ldots,\gamma_{2g}$ are the standard symplectic basis
of $H_{1}(C;\mathbb Z)$ and $\gamma_{2g+i}$ is a path connecting
$p_{1}$ and $p_{i+1}$, $i=1,\ldots, n-2g$. The period map
$\Phi:(C,\omega)\rightarrow \mathbb C^{n}$ is given by
$$ \Phi
(C,\omega)=\Big(\int_{\gamma_{1}}\omega,\ldots,\int_{\gamma_{n}}\omega\Big),$$
which provides a local coordinate system for $\mathcal
H(\mu)$, cf. \cite{Ko}. \\

Now pick a torus $T$ given by $\mathbb C/\Lambda$, where $\Lambda
= <1,\tau>$ is a lattice. Consider the coordinate $\phi=
\Phi(C,\omega)=(\phi_{1},\ldots, \phi_{n})\in \mathbb C^{n}$. We
cite the following lemma from \cite{EO}.
\begin{lemma}
$\phi_{i}\in \Lambda, i=1,\ldots, 2g$ if and only if the following conditions
hold: \\
(1) there exists a holomorphic map $f: C\rightarrow T$; \\ (2)
$\omega = f^{-1}(dz)$; \\ (3) $f$ is ramified at $q_{i}$ with
ramification order $\mu_{i}+1=l_{i}$, $i=1,\ldots, k$; \\ (4)
$f(q_{i+1})-f(q_{1})=\phi_{2g+i}$ mod $\Lambda$, $i=1,\ldots,
k-1$.
\end{lemma}

By virtue of this lemma, if we want to get a cover $C\rightarrow
T$ branched at one point with ramification type $(l_{1})\ldots
(l_{k})$, it is equivalent to $\phi_{i}=0$ mod $\Lambda$,
$i=1,\ldots, n$. Now we vary $\tau$ and such lattice points $\phi$
are dense in an open domain of $\mathbb C^{n}$. Let $U$ be the
union of such coverings $(C,\omega, q_{1},\ldots,q_{k})$, so $U$
is dense in $\mathcal H(\mu)$. Now we have the following diagram:
$$\xymatrix{
U \ar@{^{(}->}[r] \ar[drr]& \mathcal H(\mu)\ar@{^{(}->}[r] &
\mathcal
H_{g} \ar[d]^{\pi} \\
& & \mathcal M_{g} } $$ where $\mathcal H_{g}$ is the Hodge bundle
over $\mathcal M_{g}$. As long as $\mathcal H(\mu)$ dominates
$\mathcal M_{g}$, i.e., for a general $[C]\in \mathcal M_{g}$
there exists $\omega \in H^{0}(K_{C})$ such that $(\omega) =
{\displaystyle\sum_{i=1}^{k}\mu_{i}q_{i}}$, then Im
${\displaystyle\bigcup_{d}^{\infty}Y^{0}_{g,d,\sigma}} = \pi (U)$
must be dense in $\mathcal M_{g}$. Now we have to check that when
${\displaystyle\sum_{i=1}^{k}(\mu_{i}-1)\leq g-1}$, i.e., $k\geq
g-1$, $\mathcal
H(\mu)$ does dominate $\mathcal M_{g}$. \\

The case $k > g-1$ can be reduced to $k=g-1$ since there is a
natural stratification among all the moduli spaces $\mathcal
H(\mu)$. When $k=g-1$, we apply the De Jonqui\`{e}res' Formula
from \cite[VIII \S 5]{ACGH}. Suppose that $a_{1}, \ldots, a_{m}$
are distinct integers and $a_{i}$ appears $n_{i}$ times in the
partition $\mu$ of $2g-2$, ${\displaystyle \sum_{i=1}^{m}n_{i} =
g-1}$ and ${\displaystyle \sum_{i=1}^{m}n_{i}a_{i} = 2g-2}.$
Further define $R(t) =
1+{\displaystyle\sum_{i=1}^{m}a_{i}^{2}t_{i}}$ and
$P(t)=1+{\displaystyle\sum_{i=1}^{m}a_{i}t_{i}}$. Then on a genus
$g$ curve, the virtual number of canonical divisors having $n_{i}$
points of multiplicity $a_{i}$ is
$$ {\displaystyle\left[\frac{R(t)^{g}}{P(t)}\right]_{t_{1}^{n_{1}}\cdots t_{m}^{n_{m}}}}. $$ As long as this number is nonzero,
$\mathcal H(\mu)$ dominates $\mathcal M_{g}$. Let $A = {\displaystyle\sum_{i=1}^{m}a_{i}^{2}t_{i}}$ and $B = {\displaystyle\sum_{i=1}^{m}a_{i}t_{i}}$. Then
\begin{eqnarray*}
&& \left[\frac{(1+A)^{g}}{1+B}\right]_{t_{1}^{n_{1}}\cdots t_{m}^{n_{m}}} \\
& = & \left[\left(1+{g \choose 1}A+\cdots+{g\choose g-1}A^{g-1}\right)\left(1-B+B^{2}-\cdots\right)\right]_{t_{1}^{n_{1}}\cdots t_{m}^{n_{m}}} \\
& = & \left[{g\choose g-1}A^{g-1} - {g\choose g-2}A^{g-2}B + \cdots + (-B)^{g-1}\right]_{t_{1}^{n_{1}}\cdots t_{m}^{n_{m}}} \\
& = & \left[\frac{A^{g}-(A-B)^{g}}{B}\right]_{t_{1}^{n_{1}}\cdots t_{m}^{n_{m}}} \\
& = & \left[A^{g-1}+A^{g-2}(A-B) + \cdots + (A-B)^{g-1}\right]_{t_{1}^{n_{1}}\cdots t_{m}^{n_{m}}} > 0,
\end{eqnarray*}
since $A-B = {\displaystyle\sum_{i=1}^{m}(a_{i}^{2}-a_{i})s_{i}}$ has nonnegative coefficients and $A$ has positive coefficients.
\end{proof}

\begin{remark}
The most general case is $k=2g-2, \mu_{i}=1, i=1,\ldots, 2g-2$;
then we want to analyze $\alpha\beta\alpha^{-1}\beta^{-1}\in
\sigma=(2^{2g-2}1^{d-4g+4})$. Of course in this case the image of
${\displaystyle\bigcup_{d=4g-4}^{\infty}Y^{0}_{g,d,\sigma}}$ is
dense in $\mathcal M_{g}$. The most special case is $k =1,
\mu_{1}=2g-2$; then correspondingly
$\alpha\beta\alpha^{-1}\beta^{-1}\in\sigma= ((2g-1)^{1}1^{d-2g+1})$.
Obviously for dimension reasons, the image of
${\displaystyle\bigcup_{d=2g-1}^{\infty}Y^{0}_{g,d,\sigma}}$ has
to be contained in a proper subvariety of $\mathcal M_{g}$.
\end{remark}

Since the slope formula and the monodromy of $Y$ are relatively
simple to analyze when $d$ is prime, from now on we will mainly
focus on $d$ prime for concrete examples. Therefore, the following
result is also useful.

\begin{claim}
The image of ${\displaystyle\bigcup_{d\ prime}Y^{0}_{g,d,\sigma}}$
is dense in $\mathcal M_{g}$ if and only if the inequality ${\displaystyle\sum_{i=1}^{k}(\mu_{i}-1)\leq g-1}$ holds, i.e., $k\geq
g-1$.
\end{claim}

We illustrate the idea of the proof. Use the notation in the proof
of Theorem~\ref{D} and write the coordinate $\phi_{i} =
x_{i}+\sqrt{-1}y_{i}$. The degree $d$ of the map from $C$ to a
standard torus $T$ equals the area of $C$, which is given by
$$ \frac{\sqrt{-1}}{2} \int_{C}\omega \wedge \overline{\omega} =  \frac{\sqrt{-1}}{2}\sum_{i=1}^{g}(\phi_{i}\overline{\phi}_{g+i}-\overline{\phi}_{i}\phi_{g+i})
= \sum_{i=1}^{g}(x_{i}y_{g+i}-y_{i}x_{g+i}).$$ Now consider all
the integer valued vectors $(x_{1}, \ldots, x_{2g}, y_{1}, \ldots,
y_{2g})$ such that $
{\displaystyle\sum_{i=1}^{g}(x_{i}y_{g+i}-y_{i}x_{g+i})}$ is some
prime number $d$. The above claim is equivalent to the density in
$\mathbb R^{4g}$ of their union. Note for a fixed prime $d$ that
such integer points are always dense in the corresponding
hypersurface. Since there are infinitely many such hypersurfaces,
the result immediately follows.

\section{Examples of Slope Calculation}
In this section, we introduce one ad hoc method to calculate
$N_{1^{a_{1}} 2^{a_{2}}\ldots d^{a_{d}}}$ for small $g$. For
simplicity, we only deal with the case when $d$ is prime. It will
become clear that the method can also be used for
general $d$ but with more subtle analysis. \\

First, let us reduce the problem a little bit by getting rid of
the equivalence relation. As introduced before,
$Cov^{g,d,\sigma}_{1^{a_{1}} 2^{a_{2}} \dots
d^{a_{d}}}=\{(\alpha,\beta)\in S_{d}\times S_{d}\ |\
\alpha\beta\alpha^{-1}\beta^{-1}\in \sigma, \ <\alpha,\beta>
\text{is transitive}, \ \beta \in (1^{a_{1}}2^{a_{2}}\ldots d^{a_{d}})\}$.
$S_{d}$ acts on $Cov^{g,d,\sigma}_{1^{a_{1}} 2^{a_{2}} \dots
d^{a_{d}}}$ by conjugation. Burnside's lemma tells us that
$${\displaystyle N_{1^{a_{1}} 2^{a_{2}}\ldots
d^{a_{d}}}=\frac{1}{|S_{d}|}\sum_{\tau\in S_{d}}|Cov^{g,d,\sigma}_{1^{a_{1}} 2^{a_{2}} \dots d^{a_{d}}}(\tau)|},$$
where $Cov^{g,d,\sigma}_{1^{a_{1}} 2^{a_{2}} \dots d^{a_{d}}}(\tau)=\{(\alpha,\beta)\in
Cov^{g,d,\sigma}_{1^{a_{1}} 2^{a_{2}} \dots d^{a_{d}}}\ |\  \tau\alpha = \alpha\tau, \tau\beta = \beta\tau\}$.

\begin{lemma}
\label{tau1}
If $\tau$ commutes with all the elements in a transitive subgroup
$H$ of $S_{d}$, then $\tau$ must be of type $(l^{m})$ in $S_{d}$,
lm=d.
\end{lemma}

\begin{proof}
It suffices to show that for any $t\in \mathbb Z$, if $\tau^{t}$
fixes an element in $\{1,2,\ldots, d\}$, then $\tau^{t} =
id$. \\

Suppose $\tau^{t}(i)=i$. For any $j\neq i$, there exists $\xi \in H$, such
that $\xi (i)=j$. But $\tau$ also commutes with $\xi$. By
$\xi\tau^{t}\xi^{-1}=\tau^{t}$, we know that $\tau^{t}(j)=j$, so
$\tau^{t}=id$.
\end{proof}

When $d$ is prime, the above $\tau$ must be either $id$ or the
long cycle class $(d^{1})$. Pick a long cycle $(12\cdots d)$,
i.e., $i$ is sent to $i+1$. By Burnside's lemma, we have
\begin{eqnarray*}
N_{1^{a_{1}}2^{a_{2}}\ldots d^{a_{d}}} & = & \frac{1}{d!}\left(|Cov^{g,d,\sigma}_{1^{a_{1}} 2^{a_{2}} \dots d^{a_{d}}}| +(d-1)!|Cov^{g,d,\sigma}_{1^{a_{1}} 2^{a_{2}} \dots d^{a_{d}}}(12\cdots d)|\right) \\
                                                     & = & \frac{|Cov^{g,d,\sigma}_{1^{a_{1}} 2^{a_{2}} \dots d^{a_{d}}}|}{d!} + \frac{|Cov^{g,d,\sigma}_{1^{a_{1}} 2^{a_{2}} \dots d^{a_{d}}}(12\cdots d)|}{d}.
\end{eqnarray*}

\begin{lemma}
\label{tau2}
$Cov^{g,d,\sigma}_{1^{a_{1}} 2^{a_{2}} \dots d^{a_{d}}}(12\cdots d) = \emptyset$.
\end{lemma}

\begin{proof}
We know that $(12\cdots d) = \alpha (12\cdots d)\alpha^{-1} =
(\alpha (1)\alpha (2)\cdots \alpha (d))$, so $\alpha (i) = (i+s),
i = 1,2,\ldots,d$, where $s$ is an integer independent of $i$.
Similarly $\beta (j) = (j+t), j = 1,2,\ldots,d$, and $t$ is
independent of $j$. Now it is easy to check that
$\alpha\beta\alpha^{-1}\beta^{-1} = id$, so there is no solution
pair in $Cov^{g,d,\sigma}_{1^{a_{1}} 2^{a_{2}} \dots
d^{a_{d}}}(12\cdots d)$.
\end{proof}

Hence, finally we obtain the relation:
$$N_{1^{a_{1}}2^{a_{2}}\ldots d^{a_{d}}} = \frac{|Cov^{g,d,\sigma}_{1^{a_{1}} 2^{a_{2}} \dots d^{a_{d}}}|}{d!}. $$

\subsection{$g=2, \sigma = (3^{1}1^{d-3})$}
We want to find solution pairs $(\alpha, \beta)\in S_{d}\times
S_{d}$ such that $\alpha\beta\alpha^{-1}\beta^{-1}\in
(3^{1}1^{d-3})$ and $<\alpha, \beta>$ is transitive. Let $\gamma =
(abc) \in (3^{1}1^{d-3})$ satisfying the equality
$\alpha\beta\alpha^{-1}=\gamma\beta$. \\

Now the key point is, $\alpha\beta\alpha^{-1}$ and $\beta$ are in
the same conjugacy class, so $\gamma\beta \sim \beta$. Note that
$a, b, c$ cannot be contained in three different cycles of
$\beta$, since $(abc)(a\cdots) (b\cdots) (c\cdots)= (a\cdots
b\cdots c\cdots)$ changes the conjugacy type of $\beta$. So only two cases are possible: \\
(1)\ $(abc)(a\cdots b\cdots c\cdots) = (a\cdots c\cdots b\cdots)$;
\\ (2)\ $(abc)(a\cdots b\cdots )(c\cdots) = (a\cdots
c\cdots)(b\cdots)$. \\

Using also the condition that $<\alpha, \beta>$ is transitive, in
case (1) $\beta$ must be of type $(l^{m})$, so $\beta$ can only be
the long cycle $(d^{1})$. Pick an element $\beta = (12\ldots d)$.
There are ${\displaystyle{d \choose 3}}$ choices for the cycle
$(abc)$; fix one. Now for $\alpha$, we know that $(\alpha
(1)\alpha (2)\cdots \alpha (d)) = \alpha (12\cdots d) \alpha^{-1}
= (abc)(a\cdots b\cdots c\cdots) = (a\cdots c\cdots b\cdots)$, so
there are $d$ choices for $\alpha$. Overall, we get
$$N_{d^{1}} = \frac{1}{d!}(d-1)! {d \choose 3} d =
{d \choose 3}.$$

For case (2), we have $(abc)(\underbrace{a\cdots \overbrace{b\cdots}^{l_{2}}}_{l_{1}})(\underbrace{c\cdots}_{l_{2}}) = (a\cdots
c\cdots)(b\cdots)$, so $\beta$ is of type
$(l_{1}^{a_{1}}l_{2}^{a_{2}})$, $l_{1}> l_{2},
a_{1}l_{1}+a_{2}l_{2} = d$. There are $\frac{{\displaystyle
d!}}{{\displaystyle\prod_{i=1}^{2}(l_{i}^{a_{i}})(a_{i}!)}}$
choices for $\beta$. Fix one choice; then there are
$a_{1}l_{1}a_{2}l_{2}$ choices for the cycle $(abc)$. Fix one of
these also; then because of the transitivity requirement, there
are ${\displaystyle\prod_{i=1}^{2}(l_{i}^{a_{i}})(a_{i}-1)!}$
choices for $\alpha$. Multiplying all the numbers together and
dividing by $d!$, we get $$N_{l_{1}^{a_{1}}l_{2}^{a_{2}}} =
l_{1}l_{2}.$$

Therefore, we have
$$N = {d\choose 3} + \sum_{a_{1}l_{1}+a_{2}l_{2} = d \atop l_{1}>
l_{2}}l_{1}l_{2},$$ and $$ M = \frac{1}{d}{d\choose 3} +
\sum_{a_{1}l_{1}+a_{2}l_{2} = d \atop l_{1}>
l_{2}}(\frac{a_{1}}{l_{1}}+\frac{a_{2}}{l_{2}})l_{1}l_{2}.$$

\begin{theorem}
\label{g21s} When $d$ is prime and $\sigma$ is of type
$(3^{1}1^{d-3})$, the slope $s(Y_{2,d,\sigma}) = 10$.
\end{theorem}

\begin{proof}
By the slope formula, in this case
$$s(Y)=\frac{12M}{M+{\displaystyle\frac{2}{9}N}} = 10 \Longleftrightarrow \frac{N}{M} = \frac{9}{10}. $$
Actually when $d$ is prime, we have
\begin{eqnarray}
\label{Ng21}
N=\frac{3}{8}(d-2)(d-1)(d+1)
 \end{eqnarray}
 and
 \begin{eqnarray}
 \label{Mg21}
 M = \frac{5}{12}(d-2)(d-1)(d+1).
 \end{eqnarray}

We will prove these formulae in the appendix by using some tricks
from number theory.
\end{proof}

\subsection{$g=2, \sigma = (2^{2}1^{d-4})$}
We analyze this case similarly. The condition on $\alpha, \beta$
is $\alpha\beta\alpha^{-1} = (ab)(ce)\beta$. Now
there are 4 possible cases: \\
(1) $(ab)(ce)(a\cdots c\cdots b\cdots e\cdots)=(a\cdots
e\cdots b\cdots c\cdots)$; \\
(2) $(ab)(ce)(a\cdots b\cdots)(c\dots )(e\cdots) = (c\cdots e\cdots)(a\cdots )(b\cdots )$; \\
(3) $(ab)(ce)(a\cdots c\cdots )(b\cdots e\cdots) = (c\cdots
b\cdots)(e\cdots a\cdots)$;\\
(4) $(ab)(ce)(a\cdots b\cdots c\cdots)(e\cdots) =
(b\cdots e\cdots c\cdots)(a\cdots)$. \\

For case (1), $\beta$ is of type $(l^{m})$. Since $d$ is prime,
$\beta$ must be the long cycle $(d^{1})$, so there are $(d-1)!$
choices for $\beta$. Pick one; then there are  ${\displaystyle {d
\choose 4}}$ choices for the cycle $(ab)(ce)$. Fix one of these
also; then there are $d$ choices for $\alpha$. Overall, we get
$$N_{d^{1}}= {d \choose 4}.$$

For case (2), $\beta$ is of type
$(l_{1}^{a_{1}}l_{2}^{a_{2}}l_{3}^{a_{3}})$, where $l_{1}=l_{2}+l_{3}> l_{2}\geq l_{3}$. \\

(i)$\ l_{2}>l_{3}$. There are $\frac{{\displaystyle
d!}}{{\displaystyle\prod_{i=1}^{3}(l_{i}^{a_{i}})(a_{i}!)}}$
choices for $\beta$. Fix one of these; we then have
${\displaystyle\prod_{i=1}^{3}l_{i}a_{i}}$ choices for the cycle
$(ab)(ce)$. Fix one of these too. To keep the transitivity of
$<\alpha, \beta>$, there are
${\displaystyle\prod_{i=1}^{3}(l_{i}^{a_{i}})(a_{i}-1)!}$ choices
for $\alpha$. Hence, in this case we get
$$N_{l_{1}^{a_{1}}l_{2}^{a_{2}}l_{3}^{a_{3}}} = l_{1}l_{2}l_{3}.$$

(ii)$\ l_{2}=l_{3}$. Since $d$ is prime and
$a_{1}(l_{2}+l_{3})+a_{2}l_{2}+a_{3}l_{3}=d$, we get
$l_{2}=l_{3}=1$ and $l_{1}=2$. So $\beta$ is of type
$(2^{a_{2}}1^{a_{1}})$ and there are
${\displaystyle\frac{d!}{2^{a_{2}}a_{2}! a_{1}!}}$
choices. Fix $\beta$; $(ab)(ce)(ab)(c)(e)=(ce)(a)(b)$, so there
are $a_{2}{\displaystyle{a_{1}\choose 2}}$ choices for the cycle
$(ab)(ce)$. Fixing one of them, there are ${\displaystyle
2^{a_{2}+1}}(a_{2}-1)!(a_{1}-1)!$ choices for $\alpha$. Overall,
we get
$$N_{2^{a_{2}}1^{a_{1}}} = a_{1}-1.$$

For case (3), $$(ab)(ce)(\underbrace{\overbrace{a\cdots}^{p_{1}}\overbrace{ c\cdots}^{q_{1}}}_{l_{1}})(\underbrace{\overbrace{b\cdots}^{p_{2}}\overbrace{ e\cdots}^{q_{2}}}_{l_{2}}) =
(c\cdots b\cdots )(e\cdots a\cdots).$$ $l_{1}> l_{2} >1,
l_{i}=p_{i}+q_{i}$ and either $q_{1}=q_{2}$ or $p_{1}=p_{2}$.
$\beta$ is of type $(l_{1}^{a_{1}}l_{2}^{a_{2}})$: there are
${\frac{{\displaystyle
d!}}{\displaystyle\prod_{i=1}^{2}(l_{i}^{a_{i}})(a_{i}!)}}$
choices. Pick one $\beta$ and assume $p=p_{1}=p_{2}$. Since $d$ is
prime, $q_{1}\neq q_{2}$. The condition on $p$ is that $1\leq
p\leq l_{2}-1$ so there are $l_{2}-1$ choices for $p$. Fix $p$;
then there are $l_{1}l_{2}a_{1}a_{2}$ choices for the cycle
$(ab)(ce)$. Fixing one of these also, finally there are
${\displaystyle\prod_{i=1}^{2}(l_{i}^{a_{i}}) (a_{i}-1)!}$
choices for $\alpha$ . In total,
we get $l_{1}l_{2}(l_{2}-1)(d!).$ \\

For case (4), $$(ab)(ce)(\underbrace{\overbrace{a\cdots}^{l_{2}} b\cdots c\cdots}_{l_{1}} )(\underbrace{e\cdots}_{l_{2}} ) =
(a\cdots )(b\cdots e\cdots c\cdots),$$ where $l_{1}\geq l_{2}+2.$
$\beta$ has to be of type $(l_{1}^{a_{1}}l_{2}^{a_{2}})$ as in
case (3). So there are also $\frac{{\displaystyle
d!}}{{\displaystyle\prod_{i=1}^{2}(l_{i}^{a_{i}})(a_{i}!)}}$
choices for $\beta$. Fix one $\beta$; then there are
$a_{1}a_{2}l_{1}l_{2}(l_{1}-l_{2}-1)$ choices for the cycle
$(ab)(ce)$. Pick one; then there are
${\displaystyle\prod_{i=1}^{2}(l_{i}^{a_{i}})(a_{i}-1)!}$ choices
for $\alpha$. Overall, we get
$l_{1}l_{2}(l_{1}-l_{2}-1)(d!).$ \\

Now combining cases (3) and (4), finally we have
$$N_{l_{1}^{a_{1}}l_{2}^{a_{2}}} = l_{1}l_{2}(l_{1}-2).$$

Therefore, we have
$$N={d\choose 4} + \sum_{a_{1}l_{1}+a_{2}l_{2}+a_{3}l_{3}=d \atop
l_{1}=l_{2}+l_{3}>l_{2}>l_{3}} l_{1}l_{2}l_{3} $$
$$+\sum_{a_{1}l_{1}+a_{2}l_{2}=d \atop l_{1}>l_{2}}l_{1}l_{2}(l_{1}-2) +
\sum_{2a_{2}+a_{1}=d}(a_{1}-1),$$ and
$$M=\frac{1}{d}{d\choose 4} + \sum_{a_{1}l_{1}+a_{2}l_{2}+a_{3}l_{3}=d\atop
l_{1}=l_{2}+l_{3}>l_{2}>l_{3}}(\frac{a_{1}}{l_{1}}+\frac{a_{2}}{l_{2}}+\frac{a_{3}}{l_{3}})
l_{1}l_{2}l_{3} $$ $$+ \sum_{a_{1}l_{1}+a_{2}l_{2}=d \atop
l_{1}>l_{2}}(\frac{a_{1}}{l_{1}}+\frac{a_{2}}{l_{2}})l_{1}l_{2}(l_{1}-2)
+ \sum_{2a_{2}+a_{1}=d}(\frac{a_{2}}{2}+a_{1})(a_{1}-1). $$

\begin{theorem}
\label{g22s} When $d$ is prime and $\sigma$ is of type
$(2^{2}1^{d-4})$, the slope $s(Y_{2,d,\sigma}) = 10$.
\end{theorem}

\begin{proof}
By the slope formula, we know that
$$s(Y)=\frac{12M}{M+{\displaystyle\frac{N}{4}}}=10 \Longleftrightarrow
\frac{N}{M}=\frac{4}{5}.$$ Actually when $d$ is prime, we have
\begin{eqnarray}
N=\frac{1}{6}(d-3)(d-2)(d-1)(d+1)
\end{eqnarray}
and
\begin{eqnarray}
M=\frac{5}{24}(d-3)(d-2)(d-1)(d+1).
\end{eqnarray}
We also postpone the discussion of these formulae to the appendix.
\end{proof}

Combining Theorems \ref{g21s} and \ref{g22s}, we have finished the
proof of Theorem \ref{g2s} for the case $d$ prime. The general
case follows from the remark below.

\begin{remark}
\label{Mc} $\overline{\mathcal M}_{2}$ is special in that the
following equality holds:
$$\lambda = \frac{\delta_{0}}{10}+\frac{\delta_{1}}{5}.$$
Hence, for a curve $B$ in $\overline{\mathcal M}_{2}$ not entirely contained in the boundary, we always
have that its slope
$s(B)\leq10$. Moreover, $s(B) = 10$ if and only if $B \ldotp \delta_{1} = 0$. \\

There is an easy but indirect way to show that in general
$Y_{g,d,\sigma}\ldotp \delta_{i} = 0, i>0$, which was pointed out
to the author by Curtis T. McMullen. The idea is that for a cover
$\pi: C\rightarrow E$, when the vanishing cycle $\beta$ shrinks to
a node, any component $\gamma$ of $\pi^{-1}(\beta)$ is not
vanishing in $H_{1}(C)$ since $\beta$ is not a zero cycle. Then
when $\gamma$ shrinks, it will form an internal node, i.e., the
degenerate covering curve only lies in $\Delta_{0}$, but not in
$\Delta_{i}, i>0$. As in the following picture, loop (1) may
belong to $\pi^{-1}(\beta)$, but loop (2) cannot.
\begin{figure}[H]
    \centering
    \psfrag{b}{$\beta$}
    \psfrag{pi}{$\pi$}
    \psfrag{E}{$E$}
    \psfrag{C}{$C$}
    \includegraphics[scale=0.5]{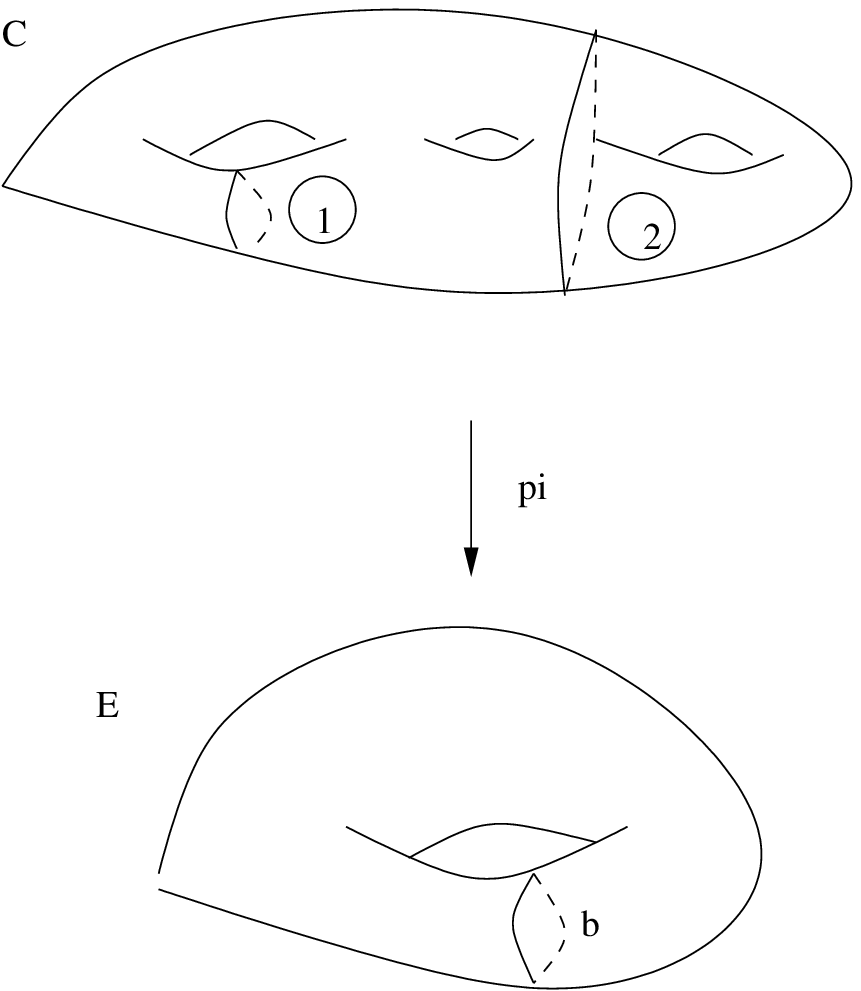}
    \label{torus}
\end{figure}

As $g=2$, $Y_{2,d,\sigma}$ may be reducible, cf. section 5,
~\cite{HL} and ~\cite{CTM1}. However, by the above argument, the
slope of each component of $Y$ is always 10.
\end{remark}

\begin{remark}
Another way to produce a one parameter family of degree $d$ covers of elliptic curves is by fixing the $j$-invariant of the target elliptic curve and moving one branched point. For instance, consider degree $d$ genus 2 covers of a fixed elliptic curve $E$
simply branched at the marked point $O$ and another point $P$. Let $P$ vary, and then we have a 1-dimensional space $W$ of admissible covers which maps to $\overline{\mathcal M}_{2}$. When $P$ meets $O$, after blowing up we will get some nodal curves
as admissible covers of $E$ with a 2-marked stable rational tail. Using the method in \cite{HM}, it is easy to write down the intersection number $W\ldotp \delta_{1}$ by the monodromy data and verify that it is not vanishing for $d\geq 3$. By the same argument in Remark \ref{Mc}, the slope of $W$ is strictly smaller than 10. Therefore, our original one parameter family $Y$ provides a better lower bound for slopes than the family $W$ does, at least for the case $g=2$.   
\end{remark}

\subsection{$g=3, \sigma = (5^{1}1^{d-5})$}
We still assume that $d$ is prime. Since the analysis is almost
the same as in the previous examples, we skip the discussion and state the result directly: \\
$$ N_{d^{1}}=8{d\choose 5}, \ N_{2^{a_{2}}1^{a_{1}}}=8(a_{2}-1),$$
$$ N_{l_{1}^{a_{1}}l_{2}^{a_{2}}} =
\frac{l_{1}l_{2}}{2}(3l_{1}^{2} +
3l_{2}^{2}-19l_{1}-11l_{2}+4d+22), l_{1}>l_{2}>1, $$
$$ N_{l_{1}^{a_{1}}l_{2}^{a_{2}}l_{3}^{a_{3}}}= \begin{cases}
11l_{1}l_{2}l_{3}, \ l_{1}\neq l_{2}+l_{3}>l_{2}>l_{3}; \\
7l_{1}l_{2}l_{3}, \ l_{1}=l_{2}+l_{3}>l_{2}>l_{3}. \end{cases}$$

Now the slope formula says that $$s(Y) =
\frac{12M}{M+{\displaystyle\frac{2N}{5}}}.$$

We calculated by computer for small prime numbers $d$ and it seems
that $s(Y)$ in this case decreases to 9. This evidence leads to
the following conjecture.

\begin{conjecture}
\label{g3s} When $\sigma$ is of type $(5^{1}1^{d-5})$, we have
$$\lim_{d\rightarrow\infty}s(Y_{3,d,\sigma})= 9.$$
\end{conjecture}

\begin{remark}
Note that for $\overline{\mathcal M}_{3}$, it is already known
that the hyperelliptic divisor $\overline H$ has the smallest
slope 9 among all effective divisors. So if the above conjecture
is true, then these curves $Y_{3,d,\sigma}$ do provide the sharp
lower bound for slopes of effective divisors on
$\overline{\mathcal M}_{3}$.
\end{remark}

Next, we study in detail the beginning case $d=5$. The result is
the following.

\begin{claim}
\label{g3d5} For $g=3, d=5$ and $\sigma$ of type $(5^{1})$,
$Y_{3,5,(5^{1})}$ has four irreducible components. Two of them have
slope $9$ and the other two have slope $9{\displaystyle
\frac{1}{3}}$.
\end{claim}

\begin{proof}
The proof is nothing but to enumerate all the possible solution
pairs $(\alpha, \beta)$ modulo equivalence relation, and further
classify the orbits by the monodromy criterion, cf. Theorem \ref{monodromy}. \\

In total there are 40 non-equivalent solutions in the following
list: \\
(1) $\alpha = (12)(34), \beta = (12345);$ (2) $\alpha = (12)(35),
\beta = (12345);$ \\
(3) $\alpha = (124), \beta = (12345);$ (4) $\alpha = (142), \beta
= (12345);$ \\
(5) $\alpha = (12453), \beta = (12345);$ (6) $\alpha = (13254),
\beta = (12345);$ \\
(7) $\alpha = (14)(25),\beta=(123);$ (8) $\alpha = (12435), \beta
= (123); $ \\
(9) $\alpha = (13425), \beta = (123);$ (10) $\alpha = (15)(23),
\beta = (12)(34);$ \\
(11) $\alpha = (135),\beta = (12)(34);$ (12) $\alpha = (12345),
\beta = (12)(34);$\\
(13) $\alpha = (12354), \beta = (12)(34);$ (14) $\alpha = (1243),
\beta = (12345);$ \\
(15) $\alpha = (1342), \beta = (12345);$ (16) $\alpha = (15)(23),
\beta = (1234);$ \\
(17) $\alpha = (15)(24), \beta = (1234);$ (18) $\alpha = (15)(34),
\beta = (1234);$ \\
(19) $\alpha = (135),\beta = (1234);$ (20) $\alpha = (125)(34),
\beta = (1234);$ \\
(21) $\alpha = (152)(34),\beta = (1234);$ (22) $\alpha = (1325),
\beta = (1234);$\\
(23) $\alpha = (1352),\beta = (1234);$ (24) $\alpha = (1523),\beta
=(1234);$ \\
(25) $\alpha = (1253),\beta = (1234);$ (26) $\alpha = (12435),
\beta = (1234);$ \\
(27) $\alpha = (14235), \beta = (1234);$ (28) $\alpha = (14)(23),
\beta = (123)(45);$ \\
(29) $\alpha = (124), \beta = (123)(45);$ (30) $\alpha = (134),
\beta = (123)(45);$ \\
(31) $\alpha = (145)(23),\beta = (123)(45);$ (32) $\alpha =
(1245), \beta = (123)(45);$ \\
(33) $\alpha = (1345), \beta = (123)(45);$ (34) $\alpha = (1425),
\beta = (123);$ \\
(35) $\alpha = (124)(35), \beta = (123);$ (36) $\alpha =
(142)(35), \beta = (123);$ \\
(37) $\alpha = (143)(25), \beta = (12)(34);$ (38) $\alpha =
(1345), \beta = (12)(34);$ \\
(39) $\alpha = (1354), \beta = (12)(34);$ (40) $\alpha = (1534),
\beta = (12)(34).$ \\

Now by monodromy action, it is routine to check that: \\
(2),(10),(13) belong to the first component $Z_{1}$; \\
(1),(3),(4),(5),(6),(7),(8),(9),(11),(12) belong to the second
component $Z_{2}$; \\
(14),(15),(16),(18),(22),(23),(24),(25),(26),(27),(38),(40) belong
to the third component $Z_{3}$; \\
(17),(19),(20),(21),(28),(29),(30),(31),(32),
(33),(34),(35),(36),(37),(39) belong to the last component
$Z_{4}$. \\

Note that the slope formula can be readily used not only for the
entire curve $Y$ but also for its irreducible components. Plugging
in the data listed above, we get
$$s(Z_{1})=s(Z_{3})=9$$ and $$ s(Z_{2})=s(Z_{4})=9{\displaystyle\frac{1}{3}}.$$
\end{proof}

For a cover $C\rightarrow E$ in $Y_{3,d,(5^{1}1^{d-5})}$, $C$
necessarily has a holomorphic 1-form with a zero of order 4. Let
$$K=\{\text{$[C]\in \mathcal M_{3}$ : $K_{C}$ has a vanishing sequence $\geq (0,1,4)$
at some point $p\in C$}\}.$$ Then we know that Im
$Y_{3,d,(5^{1}1^{d-5})}\subset K$. We define another divisor
$F\subset \mathcal M_{3}$ to be the locus of smooth plane quartics
that have hyperflexes. The slope of $\overline F$ is
$9{\displaystyle\frac{5}{8}}$, which was worked out by Cukierman,
cf. \cite{Cu}. We also know that $F\cup H = K$ and $F\cap H =
\emptyset$. It would be interesting to get some information about
the intersection of $\overline H$ and
$\overline F$ at the boundary of $\overline{\mathcal M}_{3}$. \\

First, define codimension 2 loci $$A=\{C/p\sim q: \text{$C$ is a
genus 2 curve, $p$ and $q$ are conjugate points on $C$}\},$$ and
$$B=\{\text{$E_{1}\sqcup E_{2}/p_{1}\sim p_{2}, q_{1}\sim q_{2}: E_{i}$
is an elliptic curve, and $p_{i},q_{i}\in E_{i}$}\}.$$

\begin{figure}[H]
    \centering
    \psfrag{A:}{$A:$}
    \psfrag{B:}{$B:$}
    \psfrag{E1}{$E_{1}$}
    \psfrag{E2}{$E_{2}$}
    \psfrag{C}{$C$}
    \psfrag{p~q}{$p\sim q$}
    \psfrag{p1~p2}{$p_{1}\sim p_{2}$}
    \psfrag{q1~q2}{$q_{1}\sim q_{2}$}
    \includegraphics[scale=0.6]{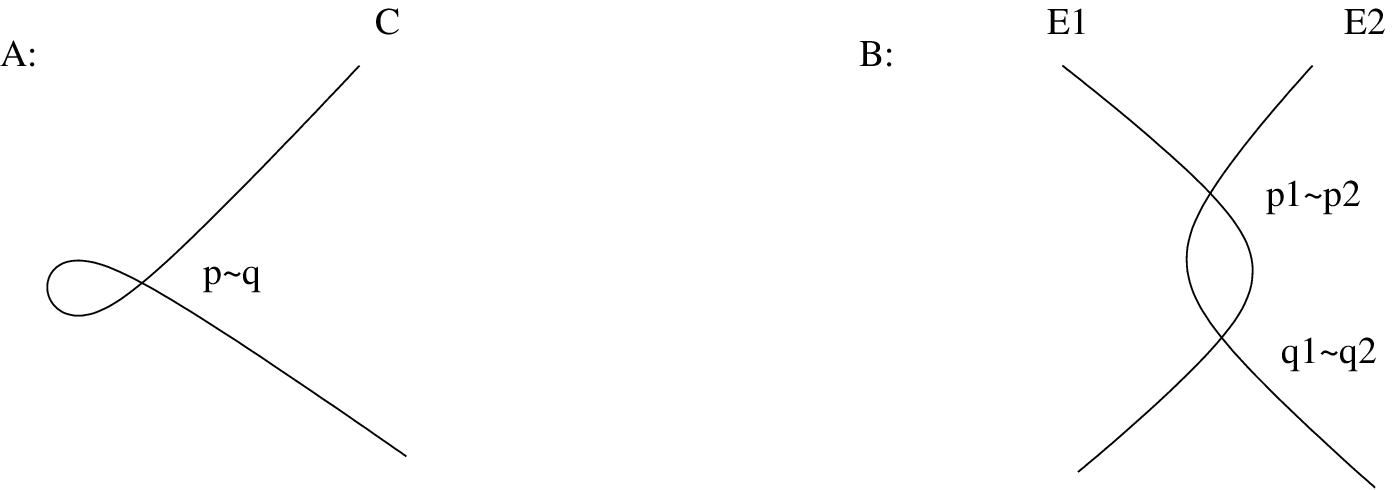}
\end{figure}

Curves in $A$ and $B$ are always double covers of rational curves
in the sense of admissible covers. It is not hard to see that
$\overline H\cap
\Delta_{0} = \overline A \cup \overline B.$ \\

For the divisor $\overline F$, we have the following claim.

\begin{claim}
$\overline F\cap \Delta_{0} \supset \overline A \cup \overline B.$
\end{claim}

\begin{proof}
First, let us verify that $A\subset F$. Take a plane conic $Q$ and
a one parameter family $C_{s}$ of general plane quartics, such
that $C_{0}$ has a hyperflex line that is also tangent to $Q$ at
the hyperflex point. By stable reduction, as in \cite[3.C]{HM1},
we know that ${\displaystyle\lim_{t\to 0}(tC_{s}+Q^{2})}$ is a
general hyperelliptic curve for $s\neq 0$. Then taking the limit
${\displaystyle\lim_{s\to 0}\lim_{t\to 0}(tC_{s}+Q^{2})}$ amounts
to squeezing two Weierstrass points together. So we get a general
element in $A$. \\

For $B$, take a banana curve $E_{1}\sqcup E_{2}/p_{1}\sim p_{2},
q_{1}\sim q_{2}$. Consider a sub linear series of $|\mathcal
O_{E_{1}}(2p_{1}+2q_{1})|$ that contracts $E_{2}$ and maps $E_{1}$
to a tacnodal plane quartic. We still have one
dimension of freedom to impose a hyperflex line. \\

Now we want to show that the inclusion in the above claim is
proper. Consider the case $g=3, d=5, \sigma = (5^{1})$ in
Claim~\ref{g3d5}. When the vanishing cycle $\beta$ is of type
$(5^{1})$, i.e., cases (1), (2), (3), (4), (5), (6), (14) and
(15), the degenerate cover is a 5-sheeted admissible covering map
from a 1-nodal geometric genus 2 curve to a rational nodal curve
totally ramified at the node and another smooth point. Note that
such a cover can be induced from a degree 5 covering map from a
smooth genus 2 curve to $\mathbb P^{1}$ totally ramified at 3
points $p, q$ and $r$. These points cannot be conjugate to each
other simultaneously. Assume that $p, q$ are not conjugate.
Identify them and identify their images also. We get a covering
curve $[C]\in$ Im $Y_{3,5,(5^{1})}\cap \Delta_{0}\subset \overline
K\cap\Delta_{0} = (\overline F\cup \overline H)\cap \Delta_{0}$.
But by the construction, we know that $[C]\notin \overline H \cap
\Delta_{0}=\overline A\cup \overline B$. Therefore, we get the
desired conclusion.
\end{proof}

\begin{remark}
It would be interesting to do the stable reduction directly for
the family $tC_{0}+Q^{2}$.
\end{remark}

\section{Counting Weighted Connected Covers}
In this section, we give a method to systematically calculate
$N_{1^{a_{1}}2^{a_{2}}\ldots d^{a_{d}}}$. An analogous exposition is givin
in the note \cite{R}. For simplicity, here we only consider the
case when $d$ is prime and $\sigma$ is of type $(2^{k}1^{d-2k}) =
\tau_{d,k}, k =
2g-2$. The general situation can be solved similarly without further difficulty. \\

For a cover $\pi: C\rightarrow E$, an automorphism $\varphi$ of
this cover is given by the following diagram:
$$\xymatrix{
C \ar[rd]_{\pi}\ar[rr]^{\varphi} & & C \ar[ld]^{\pi} \\
&  E  & }$$ If $\pi$ corresponds to one solution pair $(\alpha,
\beta) \in S_{d}\times S_{d}$, then the automorphism $\varphi$
corresponds to an element $\tau\in S_{d}$ such that
$(\tau\alpha\tau^{-1}, \tau\beta\tau^{-1}) = (\alpha, \beta)$.
Hence, we have the correspondence $Aut(C,\pi) =
Stab(\alpha,\beta)$, where $Stab(\alpha,\beta)$ is the set of stablizers of
the $S_{d}$ conjugate action. \\

In many cases people are interested in the weighted Hurwitz
numbers, i.e., counting a cover $(C,\pi)$ with weight
${\displaystyle\frac{1}{|Aut(C,\pi)|}}$. Hence, we define a
weighted number
$$\widetilde{N}^{d,k}_{\wp} =
\sum_{(\alpha,\beta)\in
Cov^{g,d,\tau_{d,k}}_{\wp}/\sim}\frac{1}{|Stab (\alpha,\beta)|}
=\frac{1}{d!}|Cov^{g,d,\tau_{d,k}}_{\wp}|, $$ where $\wp$ denotes
the conjugacy class $(1^{a_{1}}2^{a_{2}}\ldots d^{a_{d}})$. In particular
when $d$ is prime, the cover has no non-trivial automorphism since
$<\alpha,\beta>$ is transitive, cf. Lemma \ref{tau1} and \ref{tau2}. So we can pretend to
count weighted connected covers instead. \\

Now we want to get rid of the transitivity condition
imposed on $<\alpha, \beta>$. So we further define a set
$$\widehat{Cov}_{\wp}^{d,k}:= \{(\alpha,\beta)\in S_{d}\times S_{d}\ |\ \beta\in \wp, \ \alpha\beta\alpha^{-1}\beta^{-1}\in
\tau_{d,k}\},$$ and let $$\widehat{N}_{\wp}^{d,k} =
|\widehat{Cov}_{\wp}^{d,k}|.$$ We fix an element $\tau\in
\tau_{d,k}$ and look for pairs $(\gamma,\beta)\in \wp\times\wp$
such that $\gamma\beta^{-1} = \tau$. For such a pair, there are
${\displaystyle\frac{|S_{d}|}{|\wp|}}$ choices for $\alpha$ that
satisfy the equality $\alpha\beta\alpha^{-1} = \gamma$. Hence,
$$\widehat{N}_{\wp}^{d,k} =
|\tau_{d,k}|\cdot\frac{|S_{d}|}{|\wp|}\cdot
|\{(\gamma,\beta)\in \wp\times\wp: \gamma\beta^{-1}=\tau\}|.$$ By
\cite[7.68 a]{St}, we know that
$$|\{(\gamma,\beta)\in \wp\times\wp: \gamma\beta^{-1}=\tau\}| =
\frac{|\wp|^{2}}{|S_{d}|}\cdot\left(\sum_{\chi}\frac{1}{\text{deg}(\chi)}
|\chi(\wp)|^{2}\chi (\tau)\right),$$ where $\chi$ runs over all
irreducible characters of $S_{d}$. Therefore, we get the following
expression
$$\widehat{N}_{\wp}^{d,k} = |\wp|\cdot|\tau_{d,k}|\cdot \left(\sum_{\chi}\frac{1}{\text{deg}(\chi)}
|\chi(\wp)|^{2}\chi (\tau_{d,k})\right).$$

Now our task is to derive $\widetilde{N}_{\wp}^{d,k}$ from
$\widehat{N}_{\wp}^{d,k}$. Take a solution pair $(\alpha, \beta)
\in \widehat{Cov}_{\wp}^{d,k}$. The subgroup $<\alpha,\beta>$ of
$S_{d}$ may not be transitive. Consider the orbits and the action
of $\alpha,\beta$ on them, which correspond to the following data:
$$\{(\alpha_{i},\beta_{i}), i = 1,\ldots, m\ |\  
\alpha\beta\alpha^{-1}\beta^{-1}\in
(2^{k_{i}}1^{d_{i}-2k_{i}})=\tau_{d_{i},k_{i}}, $$
$$\beta_{i}\in \wp_{i}\  \text{a conjugacy class of}\  S_{d_{i}},\ 
\bigcup_{i=1}^{m}\wp_{i}=\wp,\ \sum_{i=1}^{m}k_{i} = k,$$
$$\sum_{i=1}^{m}d_{i}=d, \ <\alpha_{i},\beta_{i}>
\ \text{is a transitive subgroup of}\ S_{d_{i}} \}.$$

Two data $(\wp_{i},k_{i},d_{i})$ and $(\wp_{j},k_{j},d_{j})$ are
of the same type if $\wp_{i}\sim \wp_{j}, k_{i}=k_{j}$ and
$d_{i}=d_{j}$. Hence, we get the following equality
$$ \widehat{N}_{\wp}^{d,k} = \sum{d\choose
\underbrace{d_{1},\ldots,d_{1}}_{p_{1}},\ldots,\underbrace{d_{m},\ldots,d_{m}}_{p_{m}}}
\frac{(d_{1}!)^{p_{1}}\cdots (d_{m}!)^{p_{m}}}{(p_{1}!)\cdots
(p_{m}!)}(\widetilde{N}_{\wp_{1}}^{d_{1},k_{1}})^{p_{1}}\cdots
(\widetilde{N}_{\wp_{m}}^{d_{m},k_{m}})^{p_{m}},$$ where the
condition on the summation is
${\displaystyle\sum_{i=1}^{m}p_{i}d_{i}=d,
\sum_{i=1}^{m}p_{i}k_{i}=k}$ and
${\displaystyle\bigcup_{i=1}^{m}\wp_{i}^{p_{i}}=\wp}.$ Simplifying
the above expression, we obtain
$$ \widehat{N}_{\wp}^{d,k} = \sum (d!)
\prod_{i=1}^{m}\frac{(\widetilde{N}_{\wp_{i}}^{d_{i},k_{i}})^{p_{i}}}{p_{i}!}.$$

If $\wp$ is of type $(1^{a_{1}}\cdots d^{a_{d}})$, a datum
$(\wp_{i},k_{i},d_{i})$ corresponds to a vector
$(a_{i1},\ldots,a_{id},k_{i}), 0\leq a_{ij}\leq a_{j}, 0\leq
k_{i}\leq k, {\displaystyle \sum_{j=1}^{d}ja_{ij}=d_{i}}$ with integer entries.
Put all the vectors into a matrix $$A = \left(
\begin{array}{cccc}
a_{11} & \cdots  & a_{1d} & k_{1} \\
\vdots & \vdots & \vdots & \vdots \\
a_{m1} & \cdots  & a_{md} & k_{m}
\end{array} \right)$$
Then the summation runs over all possible $(p_{1},\ldots, p_{m})$
satisfying
$$(p_{1},\ldots,p_{m})\cdot A = (a_{1},\ldots,a_{d},k).$$

Now slightly change the notation. Define an index set $$I=
\{(a_{1},a_{2},\cdots)\ |\ a_{i}\geq 0, \ \text{and there are only
finite many non-zero entries}\}.$$ For $\wp$ of type $
(1^{a_{1}}2^{a_{2}}\cdots d^{a_{d}})$, write it as
$(1^{a_{1}}2^{a_{2}}\cdots)$. So it is determined by an element
$(a_{1},a_{2},\ldots)\in I$, and
$d=a_{1}+2a_{2}+\cdots$. \\

We define two generating functions as follows:
$$ \widehat{Z}(y;x_{1},x_{2},\ldots) =
\sum_{I,k}\frac{\widehat{N}_{\wp}^{d,k}}{(a_{1}+2a_{2}+\cdots)!}\cdot
y^{k}x_{1}^{a_{1}}x_{2}^{a_{2}}\cdots,$$ and
$$ \widetilde{Z}(y;x_{1},x_{2},\ldots) =
\sum_{I,k}\widetilde{N}_{\wp}^{d,k}\cdot
y^{k}x_{1}^{a_{1}}x_{2}^{a_{2}}\cdots,$$ where $d =
a_{1}+2a_{2}+\cdots$ and $\wp$ is of type $(1^{a_{1}}2^{a_{2}}\cdots)$ for each term. \\

Therefore, finally we obtain the relation between the generating functions for the number of connected
and possibly disconnected covers:
$$\widehat{Z} = exp(\widetilde{Z}) - 1.$$

\begin{remark}
If $k$ is odd or $2k > d = a_{1}+2a_{2}+\cdots$, then obviously
$\widehat{N}_{\wp}^{d,k} = \widetilde{N}_{\wp}^{d,k} = 0$. The
evaluation of a character $\chi$ on a conjugacy class $\wp$ can be
worked out by standard formulae from the representation theory of
$S_{d}$. However, when $d$ is large, it seems hard to evaluate the
quotient ${\displaystyle\frac{N}{M}}$ even by computer. So the
estimate of ${\displaystyle\lim_{d\to \infty}s(Y_{g,d,\sigma})}$
remains mysterious to us.
\end{remark}

\section{The Local Geometry of $Y$}
In this section, we study the local monodromy action of the map
$Y\rightarrow X$ mentioned in Remark \ref{lm}, and via it we obtain information like the
genus, orbifold points and orbifold Euler characteristic of $Y$. \\

For each component of $Y$, we take its reduced scheme structure.
An orbifold point of $Y$ possibly occurs when a smooth cover
degenerates to a singular cover, and its orbifold order depends on
the information about the degree of a local minimal base change
and the order of its extra automorphisms (not induced from nearby
covers), cf. Remark~\ref{bc}. The following two examples would
illustrate the idea.

\subsection{Basic Examples}

\textbf{Example 1: $g=2, d=3, \sigma=(3^{1})$} \\
There are three non-equivalent solution pairs $(\alpha,\beta)$: \\
(1) $\alpha = (13),\beta = (12)$; \\
(2) $\alpha = (123), \beta = (12)$; \\
(3) $\alpha = (12), \beta = (123)$. \\

When $\beta$ corresponds to the vanishing cycle, the local
monodromy action $(\alpha,\beta)\rightarrow (\alpha\beta, \beta)$ can switch sheets (1) and (2) but keep (3) unchanged,
and the global monodromy can send all three sheets to each other.
So $Y$ is a degree 3 connected cover of $X$ simply branched at
$b_{1},\ldots,b_{12}$. By Riemann-Hurwitz, $2g(Y) - 2 = 3(2g(X)-2)+ 12,$ so $g(Y) = 4.$ \\

Furthermore, over $b_{i}$, sheets (1) and (2) meet at the
ramification point $r$, and the other pre-image point $s$ lies in
sheet (3).

\begin{figure}[H]
    \centering
    \psfrag{r}{$r$}
    \psfrag{s}{$s$}
    \psfrag{bi}{$b_{i}$}
    \psfrag{X}{$X$}
    \psfrag{Y}{$Y$}
    \includegraphics[scale=0.6]{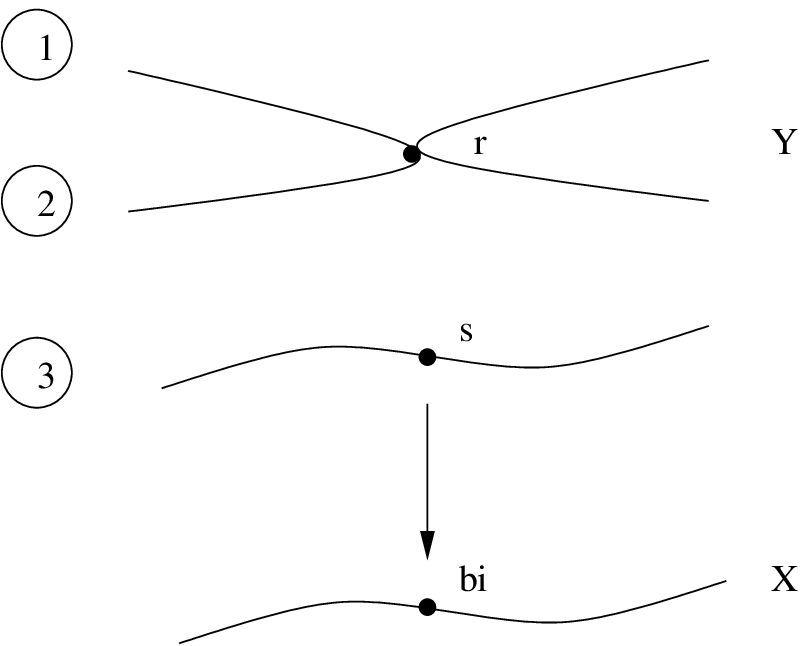}
\end{figure}

Since in case (3), $\beta$ is a length 3 cycle, as in
Remark~\ref{bc} a degree 3 base change is necessary to complete
the universal covering map locally over sheet (3). So $s$ is an
orbifold point with structure group $\mathbb Z/3$. By contrast,
locally around $r$, we need a degree 2 base change since $\beta$
contains length 2 cycles. But sheets (1) and (2) meet at $r$, so
the cover corresponding to $r$ does not have an extra order 2
automorphism compared with nearby smooth covers. Hence, $r$ is not
an orbifold point of $Y$. Finally, by the orbifold Euler
characteristic formula, we have
$$\chi(Y) = 2-2g(Y)-12(1-\frac{1}{3}) = -14.$$

\begin{remark}
By using the results about $g(Y)$ in section 5.2 and 5.3, we can
work out the orbifold Euler characteristic of $Y_{g,d,\sigma}$ in
the same way when $g=2, \sigma = (3^{1}1^{d-3})$ or
$(2^{2}1^{d-4})$ and $d$ is prime. Note that the result for
general $d$ is obtained in \cite{Ba}, although only the case
$\sigma = (3^{1}1^{d-3})$ is discussed there.
\end{remark}

It is also natural to study the degenerate covers directly and
recover the information obtained from the local monodromy action.
The cover corresponding to $s$ can be induced from a degree 3
covering map from an elliptic curve to $\mathbb P^{1}$ totally
branched over 3 points $s_{1},s_{2},p$ by identifying
$s_{1},s_{2}$ and their pre-images. The following picture reveals
the idea. A ramification point with order $k$ is marked by $(k)$
in the picture.

\begin{figure}[H]
    \centering
    \psfrag{P1}{$\mathbb P^{1}$}
    \psfrag{g=1}{$g=1$}
    \psfrag{s1}{$s_{1}$}
    \psfrag{s2}{$s_{2}$}
    \psfrag{p}{$p$}
    \psfrag{s1~s2}{$s_{1}\sim s_{2}$}
    \psfrag{3:1}{$3:1$}
    \psfrag{(3)}{$(3)$}
    \includegraphics[scale=0.6]{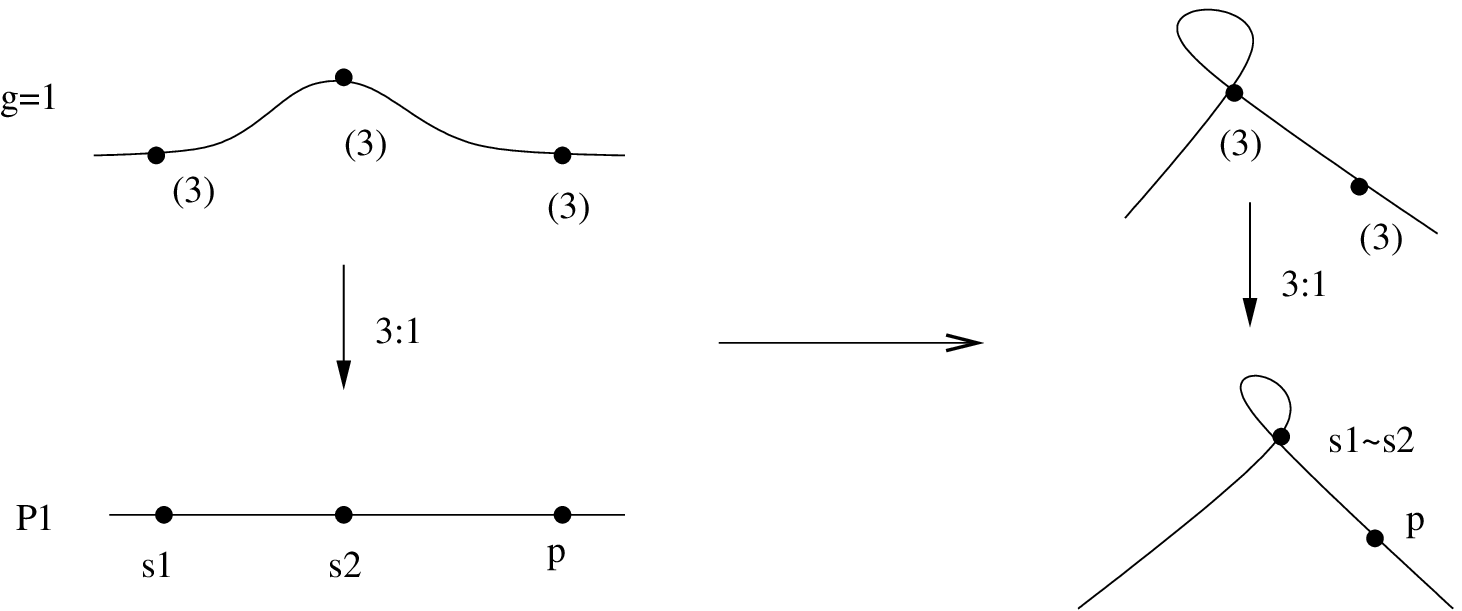}
\end{figure}

Hence, we look at triples $(\tau_{s_{1}},\tau_{s_{2}},\tau_{p})$
in $S_{3}$ such that $\tau_{s_{1}}\tau_{s_{2}}\tau_{p} = id$ and
$\tau_{i}$ is of type $(3^{1})$. The only possible solution is
$\tau_{s_{1}}=\tau_{s_{2}}=\tau_{p}$ and therefore this cover has
an order 3 automorphism, which coincides with our above discussion. \\

For $r$, we look at a degree 3 covering map $\pi$ from $\mathbb
P^{1}$ to $\mathbb P^{1}$ simply branched at $s_{1},s_{2}$ and
totally branched at $p$. Assume that
$\pi^{-1}(s_{i})=t_{i}+2t_{i}', i=1, 2$. Identify $s_{1}\sim s_{2},
t_{1}\sim t_{2}, t_{1}'\sim t_{2}'$. Then we can recover the map corresponding to $r$.

\begin{figure}[H]
    \centering
    \psfrag{t1}{$t_{1}$}
    \psfrag{t1'}{$t_{1}'$}
    \psfrag{t2}{$t_{2}$}
    \psfrag{t2'}{$t_{2}'$}
    \psfrag{s1}{$s_{1}$}
    \psfrag{s2}{$s_{2}$}
    \psfrag{p}{$p$}
    \psfrag{t1~t2}{$t_{1}\sim t_{2}$}
    \psfrag{t1'~t2'}{$t_{1}'\sim t_{2}'$}
    \psfrag{s1~s2}{$s_{1}\sim s_{2}$}
    \psfrag{P1}{$\mathbb P^{1}$}
    \psfrag{P1}{$\mathbb P^{1}$}
    \psfrag{3:1}{$3:1$}
    \psfrag{(1)}{$(1)$}
    \psfrag{(2)}{$(2)$}
    \psfrag{(3)}{$(3)$}
    \includegraphics[scale=0.6]{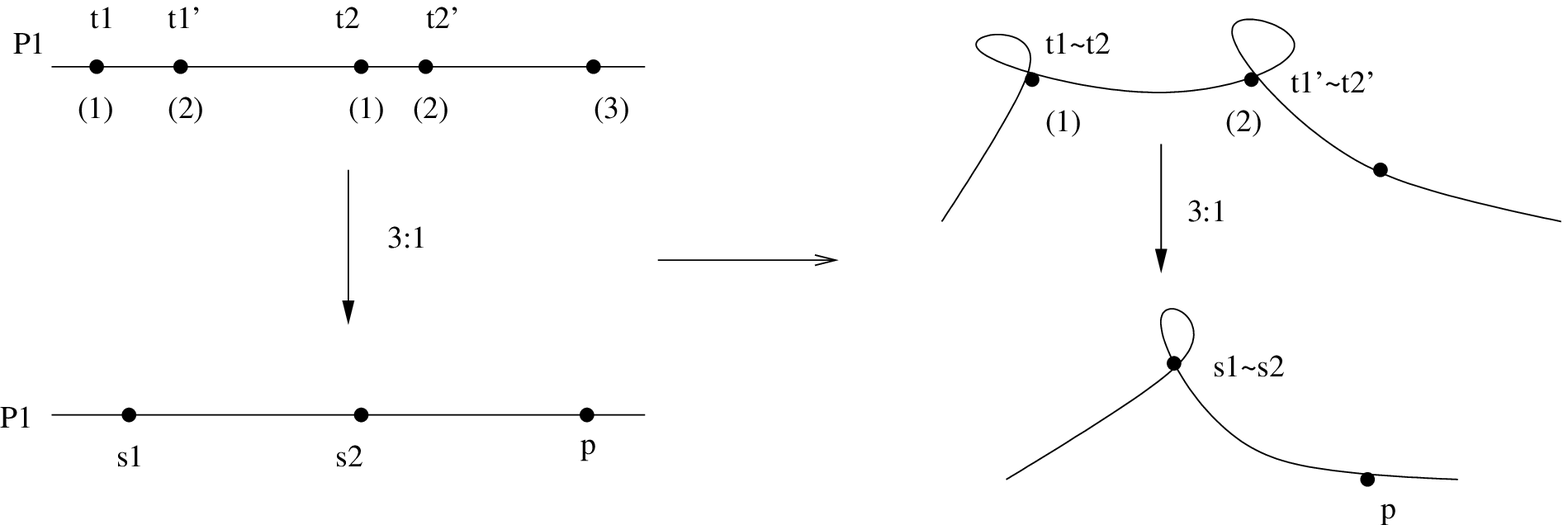}
\end{figure}

Hence, consider triples $(\tau_{s_{1}},\tau_{s_{2}},\tau_{p})$ in
$S_{3}$ such that $\tau_{s_{1}}\tau_{s_{2}}\tau_{p} = id$,
$\tau_{s_{1}}$ and $\tau_{s_{2}}$ are simple transpositions but
$\tau_{p}$ is a $(3^{1})$ cycle. Modulo the $S_{3}$ conjugation
action, there is a unique solution corresponding to
$\tau_{s_{1}}=(12),\tau_{s_{2}}=(13).$ If we switch $\tau_{s_{1}}$
and $\tau_{s_{2}}$, we get back the same cover, so it has the
automorphism induced from the involution $\iota$. But it does not
have other automorphisms. This result also coincides with our
previous analysis. \\

\textbf{Example 2: $g=3, d=5, \sigma=(5^{1})$} \\
We only focus on the case when the vanishing cycle $\beta$ is of
type $(5^{1})$, i.e., cases (1), (2), (3), (4), (5), (6), (14) and
(15) in the proof of Claim \ref{g3d5}. The local monodromy acts
transitively on (1),(3),(4),(5),(6), but keeps (2),(14),(15)
fixed. Let $r$ be the point over $b_{i}$ where the sheets
(1),(3),(4),(5),(6) meet, and $s,t,w$ be the other 3 pre-images of
$b_{i}$ contained in (2),(14),(15) respectively. Since locally
around $b_{i}$ we need a degree 5 base change, $s,t,w$ are
orbifold points with structure group $\mathbb Z/5$, and $r$ is not
an orbifold point.

\begin{figure}[H]
    \centering
    \psfrag{r}{$r$}
    \psfrag{s}{$s$}
    \psfrag{t}{$t$}
    \psfrag{w}{$w$}
    \psfrag{bi}{$b_{i}$}
    \psfrag{X}{$X$}
    \psfrag{Y}{$Y$}
    \includegraphics[scale=0.6]{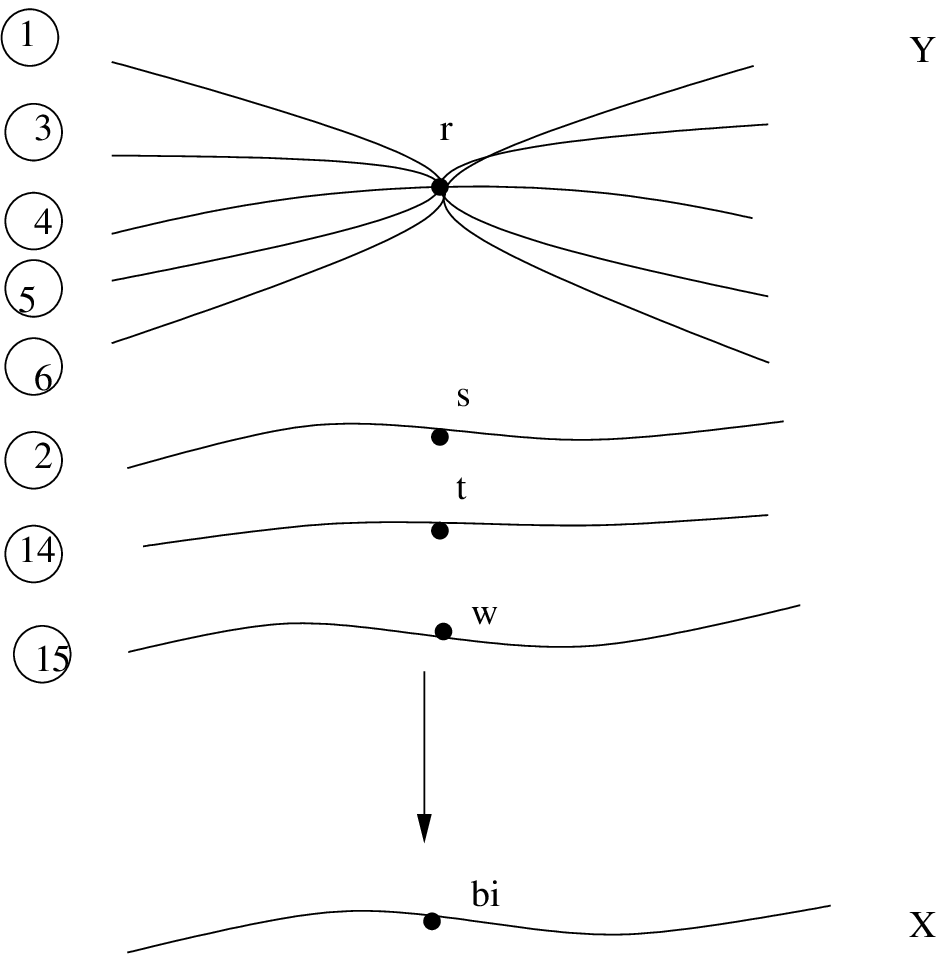}
\end{figure}

Also note that if we consider $Y'$ instead of
$Y$, then (1)-(6) are still distinct sheets in $Y'$ and they all
have automorphisms induced from the elliptic involution $\iota$.
However, (14) has to be identified with (15) under the new
equivalence relation $(\alpha,\beta)\sim (\alpha^{-1},\beta^{-1})$ for $Y'$.\\

Again, we can study the degenerate covers directly as the previous
example. In this case, the covers corresponding to $r,s,t,w$ can
be induced by degree 5 covers from genus 2 curves to $\mathbb
P^{1}$ totally branched over 3 points $s_{1},s_{2},p$. We then
identify $s_{1},s_{2}$ and their pre-image points to obtain the
desired singular covers. Now we need to analyze the solution
triples $(\tau_{s_{1}},\tau_{s_{2}},\tau_{p})$ in $S_{5}$ such that
$\tau_{s_{1}}\tau_{s_{2}}\tau_{p}=id$ and they are all $(5^{1})$
cycles. Take $\tau_{s_{1}}=(12345)$, then $\tau_{s_{2}}$ can only
be $\tau_{s_{1}}^{k},k=1,2,3$ or $(12453)$, modulo the $S_{5}$
conjugation action. For the case $\tau_{s_{2}}=\tau_{s_{1}}^{2}$,
after switching $\tau_{s_{1}}$ and $\tau_{s_{2}}$, we get a cover
equivalent to the case $\tau_{s_{2}} = \tau_{s_{1}}^{3}.$ So these two
covers can be exchanged by the involution of the target rational
nodal curve and they correspond to $t$ and $w$ in sheets (14) and
(15). For $\tau_{s_{2}}=(12453)$, switching $\tau_{s_{1}}$ and
$\tau_{s_{2}}$, we get the same cover. Moreover, this cover does
not have automorphisms except the one induced by the involution
$\iota$. So it corresponds to the non-orbifold point $r$. Finally,
if $\tau_{s_{2}} = \tau_{s_{1}}$, this cover has the automorphism
induced by $\iota$ and another order 5 automorphism, so it
corresponds to $s$ in sheet (2). \\

Now we study in general the orbits of local monodromy for the case
$g = 2$ and $d$ prime. Starting from one solution pair
$(\alpha,\beta)$ where $\beta$ is the vanishing cycle, the local
monodromy action can send $(\alpha, \beta)$ to
$(\alpha\beta^{k},\beta)$, so these two sheets have to meet at the
same degenerate cover over $b_{i}.$ We also assume that
$\alpha\beta\alpha^{-1}\beta^{-1}=\gamma\in \sigma,$ where
$\sigma$ is of type $(3^{1}1^{d-3})$ or $(2^{2}1^{d-4})$.

\subsection{$g=2, \sigma = (3^{1}1^{d-3})$}
We give a proof of Theorem~\ref{g21m}. All the numbers $N$ and
$N_{1^{a_{1}}2^{a_{2}}\ldots d^{a_{d}}}$ in the proof below are from
section 3.1.

\begin{proof}
Let $\alpha\beta\alpha^{-1}\beta^{-1}=\gamma = (abc)$ be a fixed
cycle in $S_{d}$. We look for $k$ such that there exists an
element $\tau\in S_{d}$, $\tau(\alpha\beta^{k},\beta)\tau^{-1} =
(\alpha,\beta).$ Note that such $\tau$ must satisfy
$\tau\gamma\tau^{-1}=\gamma,$ since
$$\tau\gamma\tau^{-1}=\tau(\alpha\beta^{k})\beta(\alpha\beta^{k})^{-1}\beta^{-1}\tau^{-1}
=\alpha\beta\alpha^{-1}\beta^{-1}=\gamma.$$

If $\beta$ is of type $(d^{1})$, from $\tau\beta\tau^{-1}=\beta$,
we get $\tau = \beta^{m}$ for some integer $m$. So if
$\tau(abc)\tau^{-1}=(a+m\ b+m\ c+m)=(abc)$, then $d|3m$. As long
as $d\geq 5$ is prime, $\tau$ must be $id$ and $d|k$. Since
$N_{d^{1}}={\displaystyle{d \choose 3}}$, we get
${\displaystyle\frac{1}{d}{d\choose 3}} =
{\displaystyle\frac{1}{6}}(d-1)(d-2)$ orbits each of which has
cardinality $d$. From the viewpoint of the covering map $Y\rightarrow
X$, the $d$ sheets in one orbit meet at a degenerate cover which
is not an orbifold point of $Y$, since locally around such a point
we need a degree $d$ base change to
realize the universal covering map in the proof of Theorem \ref{slope}. \\

If $\beta$ is of type $(l_{1}^{a_{1}}l_{2}^{a_{2}})$,
$l_{1}>l_{2},$ we know $(l_{1},l_{2})=1$. Without loss of
generality, assume that
$$\beta = (t_{11}t_{12}\cdots t_{1l_{1}})\cdots (t_{a_{1}1}\cdots
t_{a_{1}l_{1}})$$
$$\cdot(s_{11}s_{12}\cdots s_{1l_{2}})\cdots (s_{a_{2}1}\cdots
s_{a_{2}l_{2}}), $$ and that
$a=t_{11},b=t_{1\ l_{1}-l_{2}+1},c=s_{11}.$ From the condition
$\tau\beta\tau^{-1}=\beta,\tau\gamma\tau^{-1}=\gamma$ and $d$
prime, we can verify that $\tau$ fixes all the elements in the
cycles $(t_{11}t_{12}\cdots t_{1l_{1}})$ and $(s_{11}s_{12}\cdots s_{1l_{2}})$. 
Then we have 
$$\alpha\beta\alpha^{-1} = \gamma\beta = (t_{11}\cdots t_{1\ l_{1}-l_{2}}s_{11}\cdots s_{1l_{2}})(t_{21}\cdots t_{2l_{1}})\cdots (t_{a_{1}1}\cdots t_{a_{1}l_{1}})$$
$$\cdot (t_{1\ l_{1}-l_{2}+1}\cdots t_{1l_{1}})(s_{21}\cdots s_{2l_{2}})\cdots (s_{a_{2}1}\cdots s_{a_{2}l_{2}}).$$
So we can assume that $\alpha$ sends the cycle
$(t_{a_{1}1}\cdots t_{a_{1}l_{1}})$ to $(t_{11}\cdots
t_{1\ l_{1}-l_{2}}s_{11}\cdots s_{1l_{2}})$ and the cycle
$(s_{a_{2}1}\cdots s_{a_{2}l_{2}})$ to $(t_{1\ l_{1}-l_{2}+1}\cdots
t_{1l_{1}})$. Furthermore, assume that $\alpha (t_{ij}) = t_{i+1\ j}, 1\leq i < a_{1}-1, 1\leq j\leq l_{1}, 
\alpha (s_{ij}) = s_{i+1\ j}, 1\leq i < a_{2}-1, 1\leq j \leq l_{2},$ but $\alpha (t_{a_{1}-1\ j}) = t_{a_{1}\ j+w_{1}}$ and 
$\alpha (s_{a_{2}-1\ j}) = s_{a_{2}\ j+w_{2}}$, i.e., these actions are twisted by twist parameters $w_{1}$ and 
$w_{2}$, whose geometric meaning can be more clearly seen in \cite{HL} and the next section. 
Now, $\alpha$ contains a cycle $(t_{11}t_{21}\cdots t_{a_{1}-1\ 1}t_{a_{1}\ 1+w_{1}}\alpha(t_{a_{1}\ 1+w_{1}})\cdots )$
and the corresponding cycle in $\alpha\beta^{k}$ is $(t_{11}t_{2\ 1+k}\cdots t_{a_{1}\ 1+(a_{1}-1)k+w_{1}}\alpha(t_{a_{1}\ 1+a_{1}k+w_{1}})\cdots ).$
Since $\tau\alpha\beta^{k}\tau^{-1}=\alpha$ and $t_{11}, \alpha(t_{a_{1}\ 1+w_{1}}), \alpha(t_{a_{1}\ 1+a_{1}k+w_{1}})$ are all fixed by $\tau$, we get $l_{1}|a_{1}k$. \\

Similarly we have $l_{2}|a_{2}k$. One can check that these two conditions on $k$ are also sufficient for the existence
of $\tau$. 
Since $N_{l_{1}^{a_{1}}l_{2}^{a_{2}}}=l_{1}l_{2}$, for $\beta$ of
such type we just get $(l_{1},a_{1})(l_{2},a_{2})$ orbits and
each orbit has cardinality
${\displaystyle\frac{l_{1}l_{2}}{(l_{1},a_{1})(l_{2},a_{2})}}.$ \\

Consider the map $Y\rightarrow X$. By the Riemann-Hurwitz formula,
we have
$$2g(Y)-2=-2N+12\bigg(\frac{(d-1)(d-2)}{6}
(d-1)$$ $$+\sum_{a_{1}l_{1}+a_{2}l_{2}=d\atop
l_{1}>l_{2}}(l_{1},a_{1})(l_{2},a_{2})
\Big(\frac{l_{1}l_{2}}{(l_{1},a_{1})(l_{2},a_{2})}-1\Big)\bigg).$$ After simplifying, we get
the desired expression in the theorem. \\

Furthermore, we already know that $N\sim {\displaystyle\frac{3}{8}d^{3}}$ and
${\displaystyle\sum_{a_{1}l_{1}+a_{2}l_{2}=d\atop l_{1}>l_{2}}l_{1}l_{2}\sim
\frac{5}{24}d^{3}}$. In the appendix we will see that ${\displaystyle\sum_{a_{1}l_{1}+a_{2}l_{2}=d\atop l_{1}>l_{2}}(l_{1},a_{1})(l_{2},a_{2})}$ has lower order than $d^{3}$, so we also obtain the asymptotic result for $g(Y)$.
\end{proof}

\begin{remark}
In the next section, we will see that the above $Y$ has two
irreducible components $Z_{1}$ and $Z_{2}$ that do not intersect, cf. Theorem \ref{Curt}.
So actually we have $g(Y)=g(Z_{1})+g(Z_{2})-1.$
\end{remark}

\begin{remark}
Note that a similar genus formula for the case $g=2, \sigma =
(3^{1}1^{d-3})$ and $d$ prime is also obtained in \cite{HL} using
the technique of square-tiled surfaces. However, our space of
admissible covers $Y$ differs from the Teichm\"{u}ller discs
defined in \cite{HL} in that $Y$ is closed and it is over a pencil
of plane cubics rather than $\overline{\mathcal M}_{1,1}$.
\end{remark}

\subsection{$g=2, \sigma = (2^{2}1^{d-4})$}

Now we study the other case: $\sigma = (2^{2}1^{d-4})$, and prove
Theorem~\ref{g22m}. All the numbers $N$ and $N_{1^{a_{1}}2^{a_{2}}\ldots
d^{a_{d}}}$ in the following proof are from section 3.2.

\begin{proof}
Let $\alpha\beta\alpha^{-1}\beta^{-1}=\gamma=(ab)(ce)$ be a fixed
cycle. We still look for $k$ such that there exists $\tau\in
S_{d}$, $\tau(\alpha\beta^{k},\beta)\tau^{-1} = (\alpha,\beta).$
As before, such $\tau$ must satisfy $\tau\gamma\tau^{-1}=\gamma.$
The discussion is very similar to the one above. We analyze the
solution pairs $(\alpha,\beta)$ case by case based on the type of
the vanishing cycle
$\beta.$ \\

If $\beta$ is of type $(d^{1})$, assume that $\beta = (12\cdots
d)$. From $\tau\beta\tau^{-1}=\beta$ we get that $\tau$ must be
$\beta^{m}$ for some integer $m$. But $\tau\gamma\tau^{-1}=(a+m\
b+m)(c+m\ e+m) \neq (ab)(ce)$, since $4m \not\equiv 0\ \text{(mod
$d$)}.$ $N_{d^{1}} = {\displaystyle{d\choose 4}}$, so in this case
we get
$$\frac{1}{d}{d\choose 4} = \frac{1}{24}(d-1)(d-2)(d-3)$$
orbits, each of which contains $d$ sheets meeting at one
degenerate cover over $b_{i}$. \\

If $\beta$ is of type $(l_{1}^{a_{1}}l_{2}^{a_{2}}l_{3}^{a_{3}})$,
$l_{1}=l_{2}+l_{3}>l_{2}>l_{3},$ we can write $\beta$ as
$$(t_{11}\cdots t_{1l_{1}})\cdots (t_{a_{1}1}\cdots t_{a_{1}l_{1}})$$
$$\cdot (r_{11}\cdots r_{1l_{2}})\cdots (r_{a_{2}1}\cdots r_{a_{2}l_{2}})$$
$$\cdot (s_{11}\cdots s_{1l_{3}})\cdots (s_{a_{3}1}\cdots s_{a_{3}l_{3}}),$$
and $\gamma = (t_{11}t_{1\ l_{2}+1})(r_{11}s_{11}).$ By
$\tau\gamma\tau^{-1}=\gamma, \tau\beta\tau^{-1}=\beta$ and $l_{1},l_{2}$ co-prime, we know
that $\tau$ fixes all the elements in $(t_{11}\cdots
t_{1l_{1}}),(r_{11}\cdots r_{1l_{2}}) $ and $(s_{11}\cdots
s_{1l_{3}})$. Then we have
$$\alpha\beta\alpha^{-1} = \gamma\beta = (r_{11}\cdots r_{1l_{2}}s_{11}\cdots s_{1l_{3}})
(t_{21}\cdots t_{2l_{1}})\cdots (t_{a_{1}1}\cdots t_{a_{1}l_{1}}) $$
$$\cdot (t_{11}\cdots t_{1l_{2}})(r_{21}\cdots r_{2l_{2}})\cdots (r_{a_{2}1}\cdots r_{a_{2}l_{2}})$$
$$\cdot (t_{1\ l_{2}+1}\cdots t_{1l_{1}})(s_{21}\cdots s_{2l_{3}})\cdots (s_{a_{3}1}\cdots s_{a_{3}l_{3}}).$$
So we can assume that $\alpha$ sends the cycle
$(t_{a_{1}1}\cdots t_{a_{1}l_{1}})$ to $(r_{11}\cdots
r_{1l_{2}}s_{11}\cdots s_{1l_{3}})$, $(r_{a_{2}1}\cdots
r_{a_{2}l_{2}})$ to $(t_{11}\cdots t_{1l_{2}})$, and
$(s_{a_{3}1}\cdots s_{a_{3}l_{3}})$ to $(t_{1\ l_{2}+1}\cdots
t_{1l_{1}}).$ $\alpha(t_{ij})=t_{i+1\ j},1\leq i < a_{1}-1,1\leq
j\leq l_{1}, \alpha(r_{ij})=r_{i+1\ j},1\leq i < a_{2}-1,1\leq
j\leq l_{2}, \alpha(s_{ij})=s_{i+1\ j},1\leq i < a_{3}-1,1\leq
j\leq l_{3},$ but $\alpha(t_{a_{1}-1\ j})=t_{a_{1}\
j+w_{1}},\alpha(r_{a_{2}-1\ j})=r_{a_{2}\ j+w_{2}},
\alpha(s_{a_{3}-1\ j})=s_{a_{3}\ j+w_{3}},$ i.e., these actions
are twisted by twist parameters $w_{1},w_{2}$ and $w_{3}$.
$\alpha$ contains a cycle $(t_{11}t_{21}\cdots t_{a_{1}-1\
1}t_{a_{1}\ 1+w_{1}}\alpha (t_{a_{1}\ 1+w_{1}})\cdots )$ and the
corresponding cycle in $\alpha\beta^{k}$ is $(t_{11}t_{2\
1+k}\cdots t_{a_{1}\ 1+(a_{1}-1)k+w_{1}}\alpha(t_{a_{1}\
1+a_{1}k+w_{1}})\cdots ).$ Since
$\tau\alpha\beta^{k}\tau^{-1}=\alpha$ and
$t_{11},\alpha(t_{a_{1}1+w_{1}}),
\alpha(t_{a_{1}\ 1+a_{1}k+w_{1}})$ are all fixed by $\tau$, we get $l_{1}|a_{1}k$. \\

Similarly, we have $l_{i}|a_{i}k, i=2,3$.
These conditions are also sufficient. Since
$N_{l_{1}^{a_{1}}l_{2}^{a_{2}}l_{3}^{a_{3}}}=l_{1}l_{2}l_{3}$, we
get ${\displaystyle\prod_{i=1}^{3}(l_{i},a_{i})}$ orbits in this
case, and each orbit contains
${\displaystyle\frac{l_{1}l_{2}l_{3}}{(l_{1},a_{1})(l_{2},a_{2})(l_{3},a_{3})}}$ elements. \\

If $\beta$ is of type $(l_{1}^{a_{1}}l_{2}^{a_{2}})$,
$l_{1}>l_{2}>1$, we can write
$$\beta = (t_{11}\cdots t_{1l_{1}})\cdots (t_{a_{1}1}\cdots t_{a_{1}l_{1}})$$
$$\cdot (s_{11}\cdots s_{1l_{2}})\cdots (s_{a_{2}1}\cdots s_{a_{2}l_{2}}).$$
If $\gamma = (t_{11}s_{11})(t_{1\ p+1}s_{1\ p+1}), 0<p<l_{2}$, then
$$\alpha\beta\alpha^{-1}=\gamma\beta = (t_{1\ p+1}\cdots t_{1l_{1}}s_{11}\cdots s_{1p})(t_{21}\cdots t_{2l_{1}})
\cdots (t_{a_{1}1}\cdots t_{a_{1}l_{1}})$$
$$\cdot (t_{11}\cdots t_{1p}s_{1\ p+1}\cdots s_{1l_{2}})(s_{21}\cdots s_{2l_{2}})\cdots (s_{a_{2}1}\cdots s_{a_{2}l_{2}}).$$
We can always assume that $\alpha$ sends the cycle $(t_{a_{1}1}\cdots t_{a_{1}l_{1}})$ to
$(s_{11}\cdots s_{1p}t_{1\ p+1}\cdots t_{1l_{1}})$ and $(s_{a_{2}1}\cdots s_{a_{2}l_{2}})$ to
$(t_{11}\cdots t_{1p}s_{1\ p+1}\cdots s_{1l_{2}})$. Moreover, $\alpha(t_{ij})=t_{i+1\ j},
1\leq i <a_{1}-1$ and $\alpha(s_{ij})=s_{i+1\ j}, 1\leq i <a_{2}-1$, but
$\alpha(t_{a_{1}-1\ j})=t_{a_{1}\ j+w_{1}}$ and $\alpha (s_{a_{2}-1\ j})=s_{a_{2}\ j+w_{2}}$, where $w_{1}$ and
$w_{2}$ are the twist parameters.  \\

Now by $\tau\gamma\tau^{-1}=\gamma$, if $\tau(t_{11})=t_{1\ p+1},
\tau(s_{11})=s_{1\ p+1}, \tau(t_{1\ p+1})=t_{11}$ and $\tau(s_{1\
p+1})=s_{11}$, we get $l_{i}|2p, i=1,2$, which is impossible.
Hence, $\tau$ acts trivially on the elements in the cycles $(t_{11}\cdots t_{1l_{1}})$ and $(s_{11}\cdots s_{1l_{2}}).$ \\

The cycle of $\alpha$ starting from $t_{11}$ can be written as
$$(t_{11}t_{21}\cdots t_{a_{1}-1\ 1}t_{a_{1}\ 1+w_{1}}\alpha(t_{a_{1}\ 1+w_{1}})\cdots ),$$
and the corresponding cycle of $\alpha\beta^{k}$ is
$$(t_{11}t_{2\ 1+k}\cdots t_{a_{1}\ 1+(a_{1}-1)k+w_{1}}\alpha(t_{a_{1}\ 1+a_{1}k+w_{1}})\cdots ).$$
Since $\tau\alpha\beta^{k}\tau^{-1}=\alpha$ and $t_{11},\alpha(t_{a_{1}\ 1+w_{1}}),\alpha(t_{a_{1}\ 1+a_{1}k+w_{1}})$
are all fixed by $\tau$, we get $l_{1}|a_{1}k$. \\

Similarly, we have $l_{2}|a_{2}k$. These two conditions on $k$ are
also sufficient to find a desired $\tau$. Hence in this case,
there are $(l_{2}-1)(l_{1},a_{1})(l_{2},a_{2})$ orbits and each
orbit contains
${\displaystyle\frac{l_{1}l_{2}}{(l_{1},a_{1})(l_{2},a_{2})}}$ elements. \\

There is one more case when $\beta$ is of type
$(l_{1}^{a_{1}}l_{2}^{a_{2}})$, $l_{1}>l_{2}>1$; namely, $\gamma =
(t_{11}t_{1\ l_{2}+1})(t_{1\ l_{2}+1+m}s_{11})$ and
$$\alpha\beta\alpha^{-1} = \gamma\beta = (t_{1\ l_{2}+1}\cdots t_{1\ l_{2}+m}s_{11}\cdots s_{1l_{2}}
t_{1\ l_{2}+m+1}\cdots t_{1l_{1}})\cdots (t_{a_{1}1}\cdots t_{a_{1}l_{1}})$$
$$\cdot(t_{11}\cdots t_{1l_{2}})\cdots (s_{a_{2}1}\cdots s_{a_{2}l_{2}}).$$
We can also check that $\tau$ acts trivially on $(t_{11}\cdots t_{1l_{1}})$ and
$(s_{11}\cdots s_{1l_{2}})$. \\

Assume that, similarly to the last case, $\alpha$ acts with twist
parameters $w_{1}$ and $w_{2}$ at the end. Then $\alpha$ contains
one cycle
$$(t_{11}t_{21}\cdots t_{a_{1}-1\ 1}t_{a_{1}\ 1+w_{1}}\alpha(t_{a_{1}\ 1+w_{1}})\cdots)$$
and the corresponding cycle in $\alpha\beta^{k}$ is
$$(t_{11}t_{2\ 1+k}\cdots t_{a_{1}-1\ 1+(a_{1}-2)k}t_{a_{1}\ 1+w_{1}+(a_{1}-1)k}\alpha(t_{a_{1}\ 1+w_{1}+a_{1}k})\cdots ).$$
Since $\tau\alpha\beta^{k}\tau^{-1}=\alpha$ and
$t_{11},\alpha(t_{a_{1}\ 1+w_{1}}), \alpha(t_{a_{1}\
1+w_{1}+a_{1}k})$ are all fixed by $\tau$, we get $l_{1}|a_{1}k$.
Similarly we have $l_{2}|a_{2}k$. So in this case there are
$(l_{1}-l_{2}-1)(l_{1},a_{1})(l_{2},a_{2})$ orbits and each orbit
contains
${\displaystyle\frac{l_{1}l_{2}}{(l_{1},a_{1})(l_{2},a_{2})}}$ elements. \\

The last case is when $\beta$ is of type $(2^{a_{2}}1^{a_{1}})$.
Assume $\beta = (s_{1}t_{1})\cdots
(s_{a_{2}}t_{a_{2}})(r_{1})\cdots (r_{a_{1}})$ and $\gamma =
(s_{1}t_{1})(r_{1}r_{2}).$ Then we have
$$\alpha\beta\alpha^{-1}=\gamma\beta=(r_{1}r_{2})(s_{2}t_{2})\cdots (s_{a_{2}}t_{a_{2}})
(s_{1})(t_{1})(r_{3})\cdots (r_{a_{1}}).$$ We can further assume
$\alpha(s_{i})=s_{i+1},\alpha(t_{i})=t_{i+1}$ for $1\leq i\leq
a_{2}-1$ and $\alpha(s_{a_{2}})=r_{1},\alpha(t_{a_{2}})=r_{2}.$ Since
$\beta^{2}=id$, we only need to check when there exists $\tau$
such that $\tau\alpha\beta\tau^{-1}=\alpha$ and
$\tau\alpha\tau^{-1}=\beta.$

It is not hard to see that $\alpha$ can be of type
$$(s_{1}\cdots s_{a_{2}}r_{1}r_{3}\cdots r_{k})(t_{1}\cdots t_{a_{2}}r_{2}r_{k+1}\cdots r_{a_{1}})$$
or $$(s_{1}\cdots s_{a_{2}}r_{1}r_{3}\cdots r_{k}t_{1}\cdots
t_{a_{2}}r_{2}r_{k+1}\cdots r_{a_{1}}).$$ If $a_{2}$ is odd, then
$\alpha\beta$ can be
$$(s_{1}t_{2}\cdots s_{a_{2}}r_{2}r_{k+1}\cdots r_{a_{1}}t_{1}s_{2}\cdots t_{a_{2}}r_{1}r_{3}\cdots r_{k})$$
or $$(s_{1}t_{2}\cdots s_{a_{2}}r_{2}\cdots r_{a_{1}})(t_{1}s_{2}\cdots t_{a_{2}}r_{1}\cdots r_{k})$$ respectively.
Note that here $\alpha\beta$ is not of the same type as $\alpha$, so $\tau$ does not exist. \\

If $a_{2}$ is even, then $\alpha\beta$ can be
$$(s_{1}t_{2}\cdots s_{a_{2}-1}t_{a_{2}}r_{1}r_{3}\cdots r_{k})(t_{1}s_{2}\cdots t_{a_{2}-1}s_{a_{2}}r_{2}r_{k+1}\cdots r_{a_{1}})$$
or $$(s_{1}t_{2}\cdots s_{a_{2}-1}t_{a_{2}}r_{1}r_{3}\cdots r_{k}t_{1}s_{2}\cdots t_{a_{2}-1}s_{a_{2}}r_{2}\cdots r_{a_{1}}).$$
Here it is easy to check that the desired $\tau$ always exists. \\

Therefore, since $N_{2^{a_{2}}1^{a_{1}}}=a_{1}-1$, we get ${\displaystyle\frac{a_{1}-1}{(a_{2},2)}}$ orbits and each orbit contains
$(a_{2},2)$ elements. \\

Putting all the results together, we obtain, by the
Riemann-Hurwitz formula, for the map $Y\rightarrow X$
$$2g(Y)-2=-2N+12\bigg(\frac{1}{24}(d-1)^{2}(d-2)(d-3)+
\sum_{2a_{2}+a_{1}=d}\Big(a_{1}-1-\frac{a_{1}-1}{(a_{2},2)}\Big)
$$ $$+ \sum_{a_{1}l_{1}+a_{2}l_{2}+a_{3}l_{3}=d \atop
l_{1}=l_{2}+l_{3}>l_{2}>l_{3}}\Big(l_{1}l_{2}l_{3}-\prod_{i=1}^{3}(l_{i},a_{i})\Big)
+
\sum_{a_{1}l_{1}+a_{2}l_{2}=d\atop l_{1}>l_{2}}(l_{1}-2)\Big(l_{1}l_{2}-(l_{1},a_{1})(l_{2},a_{2})\Big)
\bigg).$$ After simplifying, we get exactly the expression in the theorem. \\

The asymptotic behavior of $g(Y)$ follows from the result about the expression of $N$.
Again, the terms involving $(l_{i},a_{i})$ do not affect the asymptotic order.
\end{proof}

\begin{remark}
In the next section we will see that the above $Y$ has at least
two components, based on the different types of subgroups
generated by $\alpha,\beta$.
\end{remark}

\section{Square-Tiled Surfaces}
In this section, we establish a correspondence between our method and the work in \cite{HL} from the viewpoint of
square-tiled surfaces. \\

The idea is quite simple. Take a standard torus $E$. If $C$ is a
cover of $E$, then $C$ can be realized as a possibly degenerate
lattice polygon with some edges and vertices identified. It covers
$d$ unit squares if the degree of the map is $d$. We will explain
the details by some examples for the case $g=2$. As before, we
only consider the situation when $d$ is prime.

\subsection{$g=2, \sigma = (3^{1}1^{d-3})$}
In this case, if there is a degree $d$ covering map $C\rightarrow
E$ only ramified at one point $q\in C$, then $C$ can be realized
as an octagon of area $d$. All of its vertices are identified to
be the unique ramification point $q$ marked with a $\bullet$ in
the following picture. Take also two loops $\alpha$ and $\beta$ of
the torus $E$ as in the picture.

\begin{figure}[H]
    \centering
    \psfrag{a}{$\alpha$}
    \psfrag{b}{$\beta$}
    \includegraphics[scale=0.7]{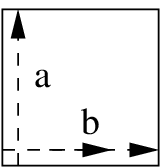}
\end{figure}

Mark the unit squares covered by the octagon by $1,2,\ldots, d$, and
consider the monodromy images in $S_{d}$ induced from $\alpha$ and $\beta$. Let us look at two examples. \\

\textbf{Example 1:} \\
Consider the following octagon.

\begin{figure}[H]
    \centering
    \psfrag{1}{$1$}
    \psfrag{2}{$2$}
    \psfrag{3}{$3$}
    \psfrag{4}{$4$}
    \psfrag{5}{$5$}
    \includegraphics[scale=0.6]{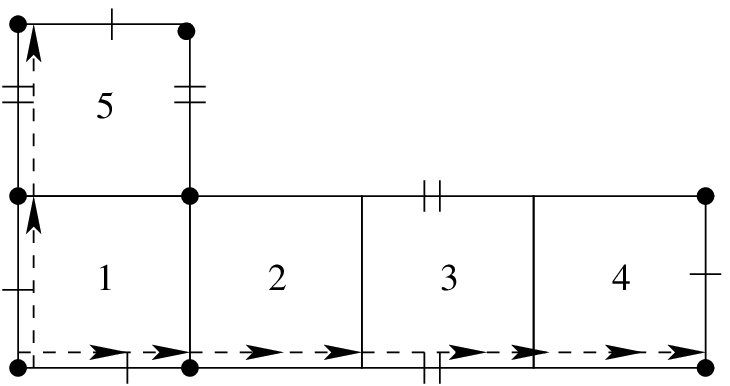}
\end{figure}

It should correspond to a degree $5$ cover of $E$. We still abuse
notation and use $\alpha, \beta$ to denote also their monodromy
images. It is easy to see $\alpha = (15)$ and $\beta = (1234).$
Then we can check that
$\alpha\beta\alpha^{-1}\beta^{-1} = (152)\in (3^{1}1^{2})$ has the desired ramification type. \\

\textbf{Example 2:} \\
Consider another octagon.

\begin{figure}[H]
    \centering
    \psfrag{1}{$1$}
    \psfrag{2}{$2$}
    \psfrag{3}{$3$}
    \psfrag{4}{$4$}
    \psfrag{5}{$5$}
    \includegraphics[scale=0.6]{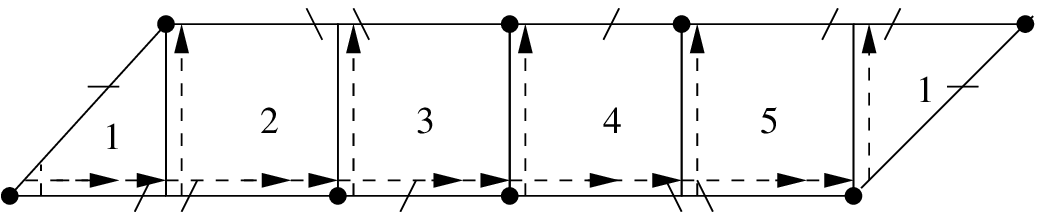}
\end{figure}

The area is still $5$. This time we get $\alpha = (12435)$ and $\beta = (12345)$.
So $\alpha\beta\alpha^{-1}\beta^{-1} = (134)\in (3^{1}1^{2})$ has the required ramification type. \\

\begin{remark}
In \cite{HL}, the square-tiled octagons can be of two types:
one-cylinder type and two-cylinder type. The one-cylinder type
corresponds to $\beta\in (d^{1})$ as in example 1, and the
two-cylinder type
corresponds to $\beta\in (l_{1}^{a_{1}}l_{2}^{a_{2}})$ as in example 2. \\

Moreover, also in \cite{HL}, it is shown that we can mark the 6
Weierstrass points of $C$. 1 or 3 out of the 6 points are integer
points, which provides two different parities invariant under the
monodromy action. It follows immediately that
$Y_{2,d,(3^{1}1^{d-3})}$ has at least two components. On the other
hand, by our method using $S_{d}$, one can check that 1 and 3
integer Weierstrass points
correspond to $<\alpha, \beta> = S_{d}$ and $A_{d}$ respectively. \\

For instance, in the first example, $<\alpha, \beta> = S_{5}$. Of
course the ramification point is one Weierstrass point, and we
mark the others with a $\Box$ in the following picture.  Note that
only the ramification point is an integer Weierstrass point.

\begin{figure}[H]
    \centering
    \includegraphics[scale=0.6]{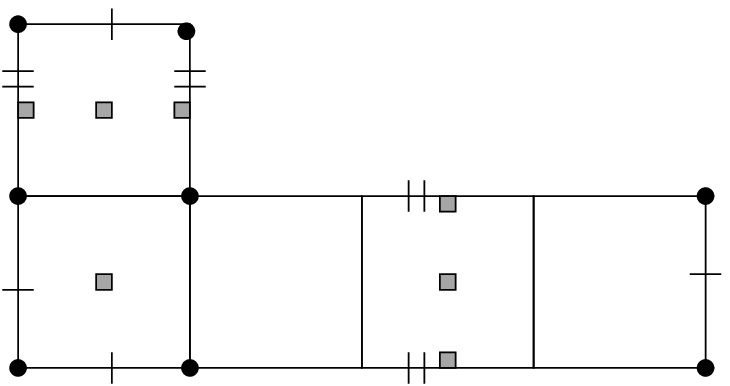}
\end{figure}

In the second example, $<\alpha, \beta> = A_{5}$. From the picture below, we can see that there are exactly 3 integer Weierstrass points.

\begin{figure}[htbp]
    \centering
    \includegraphics[scale=0.6]{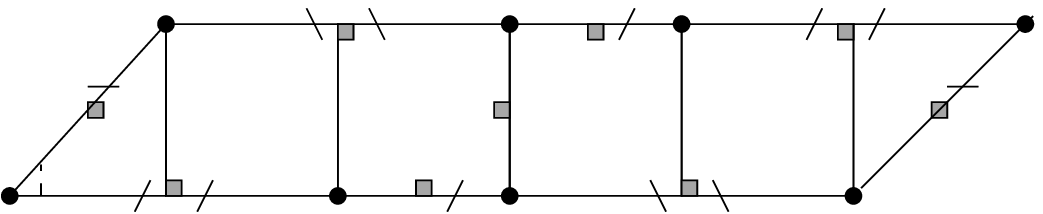}
\end{figure}

\end{remark}

\subsection{$g=2, \sigma = (2^{2}1^{d-4})$}
In this case we will have 2 ramification points on $C$. We mark
them with a $\bullet$ and a $\circ$ respectively.
As in \cite{CTM2}, $C$ can be realized as a decagon. Let us look at some examples. \\

\textbf{Example 3:} Consider a decagon in the following picture
and mark the unit squares.

\begin{figure}[H]
    \centering
    \psfrag{1}{$1$}
    \psfrag{2}{$2$}
    \psfrag{3}{$3$}
    \psfrag{4}{$4$}
    \psfrag{5}{$5$}
    \psfrag{6}{$6$}
    \psfrag{7}{$7$}
    \includegraphics[scale=0.6]{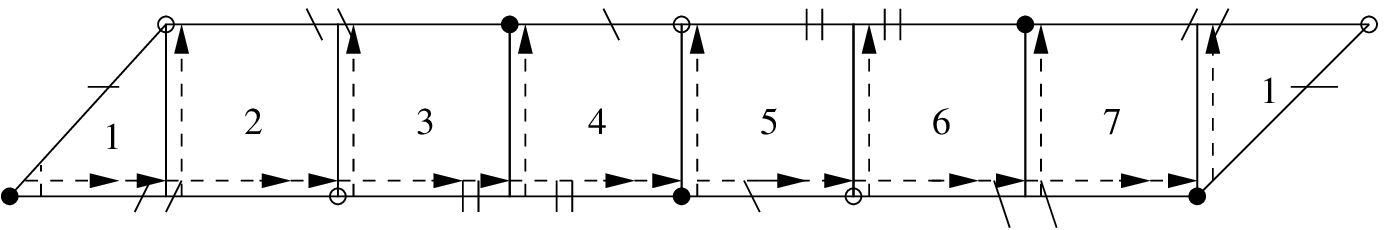}
\end{figure}

This time $d = 7.$ We get $\alpha = (1264537)$ and $\beta =
(1234567)$. So $\alpha\beta\alpha^{-1}\beta^{-1} = (16)(25)\in
(2^{2}1^{3}).$ Note that $\beta$ is of type $(d^{1})$. In general,
this is the one-cylinder type corresponding to
$\beta\in (d^{1})$ in our previous discussion. \\

\textbf{Example 4:} Consider another decagon with area 7.

\begin{figure}[H]
    \centering
    \psfrag{1}{$1$}
    \psfrag{2}{$2$}
    \psfrag{3}{$3$}
    \psfrag{4}{$4$}
    \psfrag{5}{$5$}
    \psfrag{6}{$6$}
    \psfrag{7}{$7$}
    \includegraphics[scale=0.6]{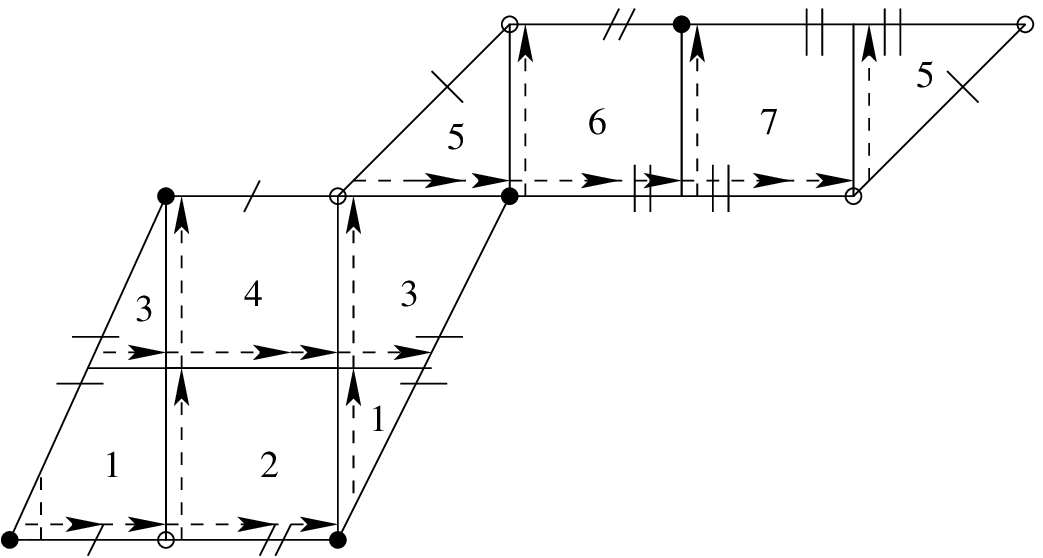}
\end{figure}

$\alpha = (1357624)$ and $\beta = (12)(34)(567)$ so
$\alpha\beta\alpha^{-1}\beta^{-1} = (16)(25)\in (2^{2}1^{3}).$ In general, this is the two-cylinder type
corresponding to $\beta\in (l_{1}^{a_{1}}l_{2}^{a_{2}})$. \\

\textbf{Example 5:} For the case $\beta\in
(l_{1}^{a_{1}}l_{2}^{a_{2}}l_{3}^{a_{3}})$, we can consider an
example such as the following.

\begin{figure}[H]
    \centering
    \psfrag{1}{$1$}
    \psfrag{2}{$2$}
    \psfrag{3}{$3$}
    \psfrag{4}{$4$}
    \psfrag{5}{$5$}
    \psfrag{6}{$6$}
    \psfrag{7}{$7$}
    \psfrag{8}{$8$}
    \psfrag{9}{$9$}
    \psfrag{10}{$10$}
    \psfrag{11}{$11$}
    \includegraphics[scale=0.6]{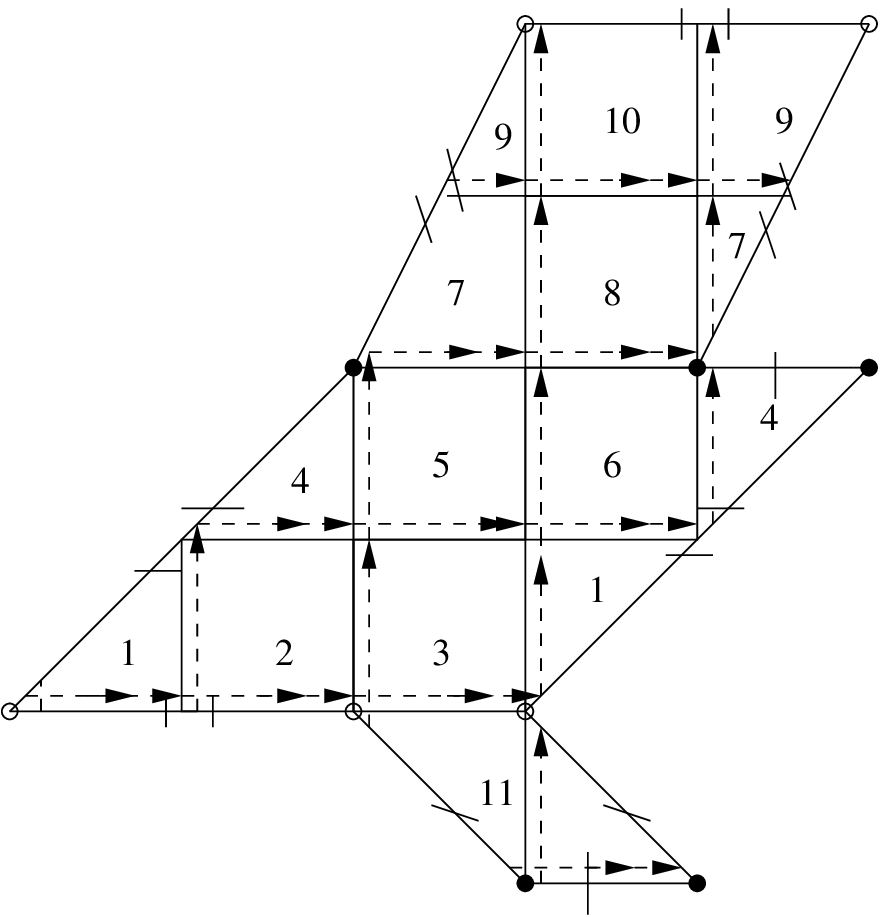}
\end{figure}

$\alpha = (168\ 10)(24\ 11\  3 5 7 9)$ and $\beta = (123)(456)(78)(9\ 10)$.
So $\alpha\beta\alpha^{-1}\beta^{-1} = (13)(7\ 11) \in (2^{2}1^{7})$. In general, this corresponds to the
three-cylinder type. \\

\begin{remark}
In \cite{HL}, there are two actions $U$ and $R$ defined by the following pictures.

\begin{figure}[H]
    \centering
    \psfrag{U}{$U$}
    \psfrag{R}{$R$}
    \psfrag{a}{$\alpha$}
    \psfrag{b}{$\beta$}
    \includegraphics[scale=0.6]{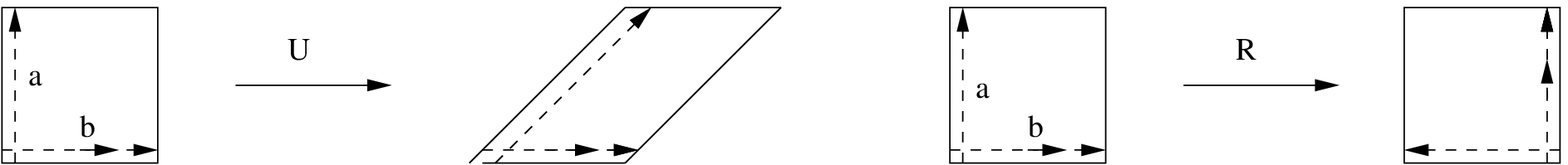}
\end{figure}

Actually they correspond to our monodromy actions $(\alpha,
\beta)\rightarrow (\alpha\beta, \beta)$ and
$(\alpha,\beta)\rightarrow (\beta^{-1},\alpha)$ respectively. For
instance, in the first example, applying the action $U$, the
direction of $\beta$ does not change but $\alpha$ changes to a
direction $\alpha'$ parallel to the diagonal.

\begin{figure}[H]
    \centering
    \psfrag{1}{$1$}
    \psfrag{2}{$2$}
    \psfrag{3}{$3$}
    \psfrag{4}{$4$}
    \psfrag{5}{$5$}
    \includegraphics[scale=0.6]{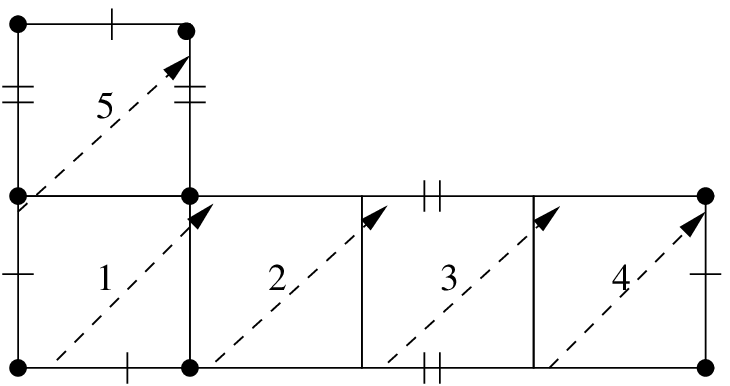}
\end{figure}

After the action $U$, we get $\alpha' = (12345) = (15)\cdot (1234)
= \alpha\beta$. Hence, the method in
\cite{HL} to work out the number of components of $Y$ by the monodromy actions can be similarly carried out here. 
\end{remark}

The reader may also be aware of the correspondence between our monodromy actions and the
butterfly moves defined in \cite{CTM1}. Actually, the result in \cite{CTM1}
is more general, not only for $d$ prime. We simply cite the result
as the following theorem.

\begin{theorem}
\label{Curt}
The Teichm\"{u}ller curve $Y_{2,d,(3^{1}1^{d-3})}$ is irreducible for $d$ even or $d=3$,
and has exactly two components for $d>3$ odd.
\end{theorem}

However, to the best of the author's knowledge, for the case $\sigma = (2^{2}1^{d-4})$ the question below is still unknown. 
\begin{question}
How many irreducible components does the curve $Y_{2,d,(2^{2}1^{d-4})}$ have? 
\end{question}

It is relatively easy to get a lower bound for the number of components. For instance, when $d > 5$ is odd, 
pick a solution pair $\alpha = (1352467\cdots d), \beta = (12)(34)$. Then $\alpha\beta\alpha^{-1}\beta^{-1} = (12)(56) 
\in \sigma = (2^{2}1^{d-4})$, and $<\alpha, \beta>$ is a subgroup of $A_{d}$. We can take another solution pair 
$\alpha = (13245\cdots d), \beta = (12)$. Then $\alpha\beta\alpha^{-1}\beta^{-1} = (12)(34)\in  \sigma = (2^{2}1^{d-4})$, and $<\alpha, \beta> = S_{d}$. Hence, in this case $Y_{2,d,(2^{2}1^{d-4})}$ has at least two components. 

\section{Appendix}
In this appendix, we will prove the equalities (1), (2), (3) and
(4) in section 3. First, we introduce the following functions
$$\sigma_{i}(n)=\sum_{k|n}k^{i}, i=1,2,\cdots.$$ These summations are quasi-modular forms with certain weights.
Now define three series
$$P = 1-24\sum_{k=1}^{\infty}\frac{kq^{k}}{1-q^{k}},$$
$$Q = 1+240\sum_{k=1}^{\infty}\frac{k^{3}q^{k}}{1-q^{k}},$$
and $$R = 1-504\sum_{k=1}^{\infty}\frac{k^{5}q^{k}}{1-q^{k}}.  $$
There are some fundamental relations among $P, Q$ and $R$ -- the
Ramanujan differential equations, cf. \cite{BY}:
$$q\frac{dP}{dq} = \frac{P^{2}-Q}{12}, \  q\frac{dQ}{dq} = \frac{PQ-R}{3}, \  q\frac{dR}{dq} = \frac{PR-Q^{2}}{2}.$$

\begin{lemma}
\label{Ra}
$$\sum_{k=1}^{d-1}\sigma_{1}(k)\sigma_{1}(d-k)=(\frac{1}{12}-\frac{d}{2})\sigma_{1}(d)+\frac{5}{12}\sigma_{3}(d).$$
\end{lemma}

\begin{proof}
Note that $$ \sum_{k=1}^{d-1}\sigma_{1}(k)\sigma_{1}(d-k) =
\left[\Big(\sum_{k=1}^{\infty}\sigma_{1}(k)q^{k}\Big)^{2}\right]_{d},$$
where $[\cdot]_{d}$ means the coefficient of the degree $d$ term
in the series expansion, and
$$\sum_{k=1}^{\infty}\sigma_{1}(k)q^{k} = \sum_{k=1}^{\infty}\sum_{j=1}^{\infty}kq^{kj} = \sum_{k=1}^{\infty}\frac{kq^{k}}{1-q^{k}}.$$
Similarly, we have
$$ \sum_{k=1}^{\infty}\sigma_{3}(k)q^{k} = \sum_{k=1}^{\infty}\frac{k^{3}q^{k}}{1-q^{k}}.$$
Hence, we get $$
\left[(\sum_{k=1}^{\infty}\sigma_{1}(k)q^{k})^{2}\right]_{d}
=\left[ \big(\frac{1-P}{24}\big)^{2}\right]_{d} =\left[
\frac{1}{24^{2}} - \frac{1}{24\cdot 12} + \frac{1}{24^{2}}(Q+12q
\frac{dP}{dq})\right]_{d}$$
$$=(\frac{1}{12}-\frac{d}{2})\sigma_{1}(d)+\frac{5}{12}\sigma_{3}(d).$$
\end{proof}

When $d$ is prime, the right side of the last equality equals
${\displaystyle\frac{1}{12}(d-1)(d+1)(5d-6)}$. Moreover,
$$\sum_{a_{1}l_{1}+a_{2}l_{2}=d, \atop l_{1}>l_{2}}l_{1}l_{2} = \frac{1}{2}\left(\sum_{k=1}^{d-1}\sigma_{1}(k)\sigma_{1}(d-k)-(d-1)\right) $$
and
\begin{eqnarray*}
\sum_{a_{1}l_{1}+a_{2}l_{2}=d, \atop l_{1}>l_{2}}(\frac{a_{1}}{l_{1}}+\frac{a_{2}}{l_{2}})l_{1}l_{2} &=& \frac{1}{2}\left(\Big(\sum_{a_{1}l_{1}+a_{2}l_{2}=d}a_{1}l_{2}+a_{2}l_{1}\Big) - d(d-1)\right) \\
 &=& \bigg(\sum_{a_{1}l_{1}+a_{2}l_{2}=d}l_{1}l_{2}\bigg) - \frac{1}{2}d(d-1).
\end{eqnarray*}

Now the equalities (1), (2) for $N$ and $M$ for $g=2, \sigma$ of
type $(3^{1}1^{d-3})$ follow immediately.
Similarly one can verify the equalities (3) and (4). \\


Now we turn to the asymptotic result for $g(Y)$ in
Theorem~\ref{g21m}. It suffices to verify that ${\displaystyle
\sum_{a_{1}l_{1}+a_{2}l_{2}=d}}(l_{1}, a_{1})(l_{2}, a_{2})$ has
asymptotic order less than 3. Actually in \cite{HL}, the order is
proved to be less than ${\displaystyle\frac{3}{2}}+\varepsilon$
for any $\varepsilon > 0.$ The asymptotic behavior of $g(Y)$ in
Theorem~\ref{g22m} can be estimated similarly to show that it is
less than 4.

Department of Mathematics, Harvard University, 1 Oxford Street, Cambridge, MA 02138 \par
{\it Email address:} dchen@math.harvard.edu

\end{document}